\documentclass[12pt]{degm}

\usepackage[latin1]{inputenc}
\usepackage{amsfonts}
\usepackage{amsmath}
\usepackage{amssymb}
\usepackage[all]{xy}

\newtheorem{theo}{Th{\'e}or{\`e}me}[section]
\newtheorem{prop}[theo]{Proposition}
\newtheorem{cor}[theo]{Corollaire}
\newtheorem{lemme}[theo]{Lemme}

\newtheorem{rque}[theo]{Remarque}
\newtheorem{defin}[theo]{D{\'e}finition}
\newtheorem{ex}[theo]{Exemple}

\renewcommand\refname{R{\'e}f{\'e}rences}

\DeclareMathOperator{\bs}{\backslash}

\DeclareMathOperator{\an}{an}

\DeclareMathOperator{\Aff}{Aff}
\DeclareMathOperator{\Alg}{Alg}

\DeclareMathOperator{\ev}{ev}
\DeclareMathOperator{\End}{End}

\DeclareMathOperator{\Form}{Form}

\DeclareMathOperator{\gal}{gal}

\DeclareMathOperator{\GL}{GL}
\DeclareMathOperator{\G}{G}
\DeclareMathOperator{\Hom}{Hom}

\DeclareMathOperator{\im}{im}

\DeclareMathOperator{\Mod}{mod}

\DeclareMathOperator{\N}{N}

\DeclareMathOperator{\Proj}{Proj}

\DeclareMathOperator{\Res}{Res}
\DeclareMathOperator{\Rig}{Rig}
\DeclareMathOperator{\rig}{rig}

\DeclareMathOperator{\SL}{SL}
\DeclareMathOperator{\Spec}{Spec}
\DeclareMathOperator{\Spf}{Spf}
\DeclareMathOperator{\Sprig}{Sp_{rig}}

\DeclareMathOperator{\Sym}{Sym}
\DeclareMathOperator{\Tr}{Tr}

\DeclareMathOperator{\p}{\mathbb{P}}

\DeclareMathOperator{\Z}{\mathbb{Z}}
\DeclareMathOperator{\n}{\mathbb{N}}
\DeclareMathOperator{\Gm}{\mathbb{G}_m}
\DeclareMathOperator{\Q}{\mathbb{Q}}
\DeclareMathOperator{\R}{\mathbb{R}}
\DeclareMathOperator{\C}{\mathbb{C}}

\title{\Large Compactifications arithm{\'e}tiques des vari{\'e}t{\'e}s de Hilbert et 
formes modulaires de Hilbert pour $\Gamma_1(\mathfrak{c},\mathfrak{n})$} 
 
\author{\large Mladen Dimitrov}
\date{\today}

\begin{document}

\maketitle

Soit $F$   un   corps  de nombres totalement r{\'e}el de degr{\'e} $d_F$, 
d'anneau des entiers $\mathfrak{o}$, de diff{\'e}rente $\mathfrak{d}$
et de discriminant $\Delta_F=\N_{F\!/\!\Q}(\mathfrak{d})$.
 On abr{\'e}gera $\N=\N_{F\!/\!\Q}$.

On se donne un groupe alg{\'e}brique $D_{/\Q}$, interm{\'e}diaire entre $\Gm$ et 
$\Res^F_{\Q} \Gm$, connexe :
$ \Gm\hookrightarrow D\hookrightarrow \Res^F_{\Q} \Gm$.
On d{\'e}finit le groupe alg{\'e}brique  $G^D_{/\Q}$ 
 (resp. $G^*_{/\Q}$) comme le  produit fibr{\'e} de 
$D$ (resp. $\Gm$) et de  $\Res^F_{\Q} \GL_2$ au-dessus de  $\Res^F_{\Q} \Gm$.
On  a le  diagramme cart{\'e}sien suivant :

\vspace{-3mm}
$$\xymatrix@R=10pt{\Res^F_{\Q} \SL_2 \ar[d]\ar@{^{(}->}[r] & G^* \ar[d]\ar@{^{(}->}[r] & 
G^D \ar[d]\ar@{^{(}->}[r] & \Res^F_{\Q} \GL_2\ar[d]^{\nu} \\
1 \ar@{^{(}->}[r]& \Gm \ar@{^{(}->}[r]& D \ar@{^{(}->}[r]& \Res^F_{\Q} \Gm,}$$

\vspace{-3mm}
\noindent o{\`u} la fl{\`e}che $\nu:\Res^F_{\Q} \GL_2\rightarrow \Res^F_{\Q} \Gm$ 
est  donn{\'e}e par la norme r{\'e}duite.

\smallskip

Le sous-groupe de Borel standard de $G^D$, 
son radical unipotent et  son tore maximal standard sont not{\'e}s 
$B$, $U$ et  $T$, respectivement. On pose $T_1=T\cap \ker(\nu)$. 

Pour toute $\Q$-alg{\`e}bre $R$ et pour tout groupe alg{\'e}brique $H$ sur $\Q$, on 
note $H_R$ le groupe de ses $R$-points. 

Soit $\mathfrak{n}$  un id{\'e}al de $\mathfrak{o}$ premier {\`a} $\Delta_F$
et ne divisant ni 2, ni 3 et soit   $\mathfrak{c}$   un id{\'e}al fractionnaire 
de $F$, que l'on peut supposer premier {\`a} $\mathfrak{n}$.   
Alors le groupe de congruences 
$\Gamma=\Gamma_1^D(\mathfrak{c},\mathfrak{n})$, d{\'e}fini dans la
partie 3, est sans torsion et l'espace 
 de modules de vari{\'e}t{\'e}s ab{\'e}liennes de Hilbert-Blumenthal correspondant 
  $M=M_1^D(\mathfrak{c},\mathfrak{n})$  est un 
$\Z[\frac{1}{\N(\mathfrak{n})}]$-sch{\'e}ma,  lisse au-dessus 
de  $\Z[\frac{1}{\Delta}]$, o{\`u}
$\Delta=\N(\mathfrak{dn})$ (voir la partie \ref{carte-locale} pour une d{\'e}finition 
pr{\'e}cise de l'espace de modules $M$).

\medskip
Cet article d{\'e}crit les compactifications  arithm{\'e}tiques de $M$ 
et donne quelques unes de leurs  propri{\'e}t{\'e}s.

Les principales r{\'e}f{\'e}rences sont  les articles  \cite{rapoport} de M. Rapoport
et \cite{chai} de C.-L. Chai, o{\`u} les compactifications toro{\"\i}dales et minimale 
sont construites pour le sous-groupe de congruence principal de 
niveau $\N(\mathfrak{n})$, lorsque $D=\Gm$.  Par ailleurs,  Rapoport
explique comment on peut obtenir une compactification partielle de 
$M$ aux pointes non-ramifi{\'e}es. 
La  contribution principale de ce travail est qu'il  fournit les
cartes locales servant {\`a} compactifier les pointes ramifi{\'e}es. 
Une  application imm{\'e}diate  est le ``principe du $q$-d{\'e}veloppement''
en ces pointes ramifi{\'e}es.

Les r{\'e}sultats de cet article sont utilis{\'e}s dans un article   commun 
avec J. Tilouine \cite{dimtildg}, o{\`u} figurent aussi  diff{\'e}rentes   
applications   aux formes modulaires de Hilbert. En vue  de ces 
applications,  il est important 
de disposer de  compactifications toro{\"\i}dales lisses  
de $M$, puisque  l'on sait prolonger les fibr{\'e}s automorphes {\`a} celles-ci.

Le groupe auxiliaire $D$ nous permet  de traiter simultan{\'e}ment 
le cas du groupe modulaire de Hilbert et celui de sa version {\'e}tendue, qui 
sont d'{\'e}gale importance et correspondent 
 {\`a} $D=\Gm$ et $D=\Res^F_{\Q} \Gm$, respectivement (voir $\cite{BL}$). 

\smallskip
Je remercie tous ceux  qui m'ont consacr{\'e} du temps pour  discuter, 
et en particulier Y. Henrio, qui a eu la gentillesse de m'expliquer le
th{\'e}or{\`e}me de descente formelle de Rapoport, ainsi que 
A. Abb{\`e}s, D. Barsky, G. Chenevier, H. Hida, A. Mokrane, M. Raynaud
 et  E. Urban. 
Je voudrais  exprimer toute ma gratitude {\`a} J. Tilouine parce qu'il  
m'a initi{\'e} {\`a}  ce sujet de recherche passionnant et  constamment 
encourag{\'e} au cours de la pr{\'e}paration de ce travail.
Enfin, je remercie  les rapporteurs pour leurs remarques int{\'e}ressantes. 

\smallskip
Nous rappelons d'abord bri{\`e}vement la construction g{\'e}n{\'e}rale de 
vari{\'e}t{\'e}s semi-ab{\'e}liennes, donn{\'e}e par D. Mumford dans 
le cas totalement d{\'e}g{\'e}n{\'e}r{\'e} \cite{mum}. Nous introduisons ensuite 
la notion de $(R,\mathfrak{n})$-pointe, version alg{\'e}brique 
de la $\Gamma$-pointe. Cela nous permet 
de construire, en suivant \cite{rapoport}, les
cartes locales, qui seront utilis{\'e}es pour
les compactifications toro{\"\i}dales arithm{\'e}tiques.

\section{La construction de Mumford.}

Soit $R$ un anneau excellent, int{\'e}gralement clos, noeth{\'e}rien, complet 
pour la  topologie $I$-adique, pour un   id{\'e}al radiciel $I=\sqrt{I}$.
Soit $K$ le corps des fractions de $R$. 

Soit $S=\Spec(R)$, $\eta$ son point g{\'e}n{\'e}rique  et 
 $S_0=\Spec(R/I)$ le sous-sch{\'e}ma ferm{\'e} d{\'e}fini par $I$.

\begin{defin}
Un $S$-sch{\'e}ma en groupes commutatif, lisse et de type fini $G$ est dit 
semi-ab{\'e}lien,  si ses fibres g{\'e}om{\'e}triques sont des 
 extensions d'une vari{\'e}t{\'e} ab{\'e}lienne par un tore. 
\end{defin}

Consid{\'e}rons le tore d{\'e}ploy{\'e} $\widetilde{G}=\mathbb{G}_m^r\times S$  
de rang $r$ sur $S$. Soit $\mathfrak{b}$ un sous-groupe discret 
polarisable de $\widetilde{G}_\eta$.
L'objet de cette section est d'esquisser la construction 
d'un sch{\'e}ma semi-ab{\'e}lien $G/S$, comme ``quotient'' de $\widetilde{G}$
par  $\mathfrak{b}$. La strat{\'e}gie est la suivante :

\smallskip

(i) Construire une ``{compactification}'' $\widetilde{G}\hookrightarrow \widetilde{P}$
telle que l'action de $\mathfrak{b}$ s'{\'e}tende {\`a} $\widetilde{P}$ et que $\mathfrak{b}$ agisse 
librement et discontinument sur 
$\widetilde{P} \underset{S}{\times}S_0$ (pour la topologie de Zariski).

(ii) Suivre les fl{\`e}ches du diagramme :
\xymatrix@R=15pt{\widetilde{G}\ar@<-2pt>@{^{(}->}^{\text{ouvert}}[rr]& &  \widetilde{P} & & 
\widetilde{\mathfrak{P}}.
\ar@<1pt>_{\text{compl{\'e}tion}}[ll]
\ar_{\text{quotient formel par } \mathfrak{b}}[d]\\
G\ar@<-1pt>@{^{(}->}^{\text{ouvert}}[rr]&&P
 && \mathfrak{P}\ar@<1pt>^{\text{alg{\'e}brisation}}[ll]  }

(iii) Enfin, montrer que $G$ est semi-ab{\'e}lien sur $S$, ind{\'e}pendant 
du choix de $\widetilde{P}$, que $G_\eta$ est ab{\'e}lien,  et que 
$G_0=\widetilde{G}_0=\mathbb{G}_m^r\times S_0$.

\smallskip

\noindent{\bf P{\'e}riodes et polarisation.}
Soit $\mathfrak{a}= \Z^r$ le
groupe des caract{\`e}res de $\widetilde{G}$.
Pour $\alpha\in \mathfrak{a}$, notons ${\mathfrak{X}}^\alpha\in 
\mathrm{H}^0(\widetilde{G},  \mathcal{O}_{\widetilde{G}})$
le caract{\`e}re associ{\'e}. Alors de mani{\`e}re canonique :

\vspace{-3mm}
$$\widetilde{G}=\Spec(R[\mathfrak{X}^\alpha;\alpha\in \mathfrak{a}])$$
\vspace{-2mm}
\begin{defin}
Un ensemble de {\it p{\'e}riodes} est un sous-groupe $\mathfrak{b} 
\subset  \widetilde{G}_\eta$ isomorphe {\`a} $\Z^r$.
\end{defin}

\vspace{-2mm}
\begin{defin}\label{polar}
Une {\it polarisation} pour $\mathfrak{b}$ est un homomorphisme 
$\phi:\mathfrak{b}\rightarrow \mathfrak{a}$
tel que :

$\mathrm{(i)}$  
$\mathfrak{X}^{\phi(\beta)}(\beta')=\mathfrak{X}^{\phi(\beta')}(\beta)$, pour tout  
$\beta,\beta'\in \mathfrak{b}$,

$\mathrm{(ii)}$ 
$\mathfrak{X}^{\phi(\beta)}(\beta)\in I$, pour tout  $\beta\in \mathfrak{b}\bs \{0\}$.
\end{defin}
\begin{lemme}\label{lemme2}
Pour tout $\alpha\in \mathfrak{a}$, il existe un entier $n\geq 1$ avec $\mathfrak{X}^{n\phi(\beta)+\alpha}(\beta)\in R$
pour tout $\beta\in\mathfrak{b}$.
\end{lemme}

\noindent{\bf Mod{\`e}les relativement complets.}
{\'E}tant donn{\'e} un ensemble de p{\'e}riodes $\mathfrak{b}
\subset  \widetilde{G}_\eta$ muni d'une polarisation $\phi$,
Mumford donne la

\vspace{-3mm}
\begin{defin}
Un {\it mod{\`e}le relativement complet} de $\widetilde{G}$, par rapport {\`a} $(\mathfrak{b},\phi)$, 
est la donn{\'e} des {\'e}l{\'e}ments suivants :

\noindent$\mathrm{(a)}$  Un sch{\'e}ma int{\`e}gre $\widetilde{P}$, localement de type fini sur $R$,

\noindent$\mathrm{(b)}$  Une immersion ouverte $i:\widetilde{G}\hookrightarrow\widetilde{P}$,

\noindent$\mathrm{(c)}$  Un faisceau inversible $\widetilde{\mathcal{L}}$ sur $\widetilde{P}$,

\noindent$\mathrm{(d)}$  Une action du tore 
$\widetilde{G}$ sur $\widetilde{P}$ et $\widetilde{\mathcal{L}}$,
not{\'e}e $S_g:\widetilde{P}\rightarrow\widetilde{P}$ et 
$S_g^*:\widetilde{\mathcal{L}}\rightarrow\widetilde{\mathcal{L}}$, pour tout 
point fonctoriel $g$ de $\widetilde{G}$,

\noindent$\mathrm{(e)}$  Une action de $\mathfrak{b}$ sur $\widetilde{P}$ et $\widetilde{\mathcal{L}}$,
not{\'e}e $T_{\beta}:\widetilde{P}\rightarrow\widetilde{P}$ et 
$T_{\beta}^*:\widetilde{\mathcal{L}}\rightarrow\widetilde{\mathcal{L}}$, pour
tout $\beta\in \mathfrak{b}$,

\smallskip
  satisfaisant aux conditions suivantes :
\smallskip

\noindent$\mathrm{(i)}$  Il existe un ouvert $\widetilde{G}$-invariant $U\subset \widetilde{P}$
de type fini sur $S$ et tel que $\widetilde{P}=\cup_{\beta\in\mathfrak{b}}T_{\beta}(U)$.

\noindent$\mathrm{(ii)}$  Pour toute valuation $v$ sur le corps des
fonctions rationnelles sur $\widetilde{G}$ et  qui est positive sur $R$, on a :

\noindent $v$ a du centre sur $\widetilde{P}$  $\iff$ pour tout 
$\alpha\in \mathfrak{a}$, il existe  $\beta\in\mathfrak{b}$ avec 
$v(\mathfrak{X}^{\alpha}(\beta)\mathfrak{X}^{\alpha})\geq 0$.

\noindent$\mathrm{(iii)}$  Les actions de $\widetilde{G}$ et  $\mathfrak{b}$ sur $\widetilde{P}$ 
prolongent leurs actions par translation sur $\widetilde{G}_{\eta}$.

\noindent$\mathrm{(iv)}$   Les actions de $\widetilde{G}$ et  $\mathfrak{b}$ sur $\widetilde{\mathcal{L}}$ 
v{\'e}rifient la condition de compatibilit{\'e} suivante :

$S_g^*T_{\beta}^*=\mathfrak{X}^{\phi(\beta)}(g)T_{\beta}^*S_g^*$, pour tout $\beta\in\mathfrak{b}$ et tout point fonctoriel  $g$ de $\widetilde{G}$.

\noindent$\mathrm{(v)}$  $\widetilde{\mathcal{L}}$ est ample sur  $\widetilde{P}$, au sens 
que les compl{\'e}ments des lieux des z{\'e}ros des sections globales 
 $ \mathrm{H}^0(\widetilde{P},\widetilde{\mathcal{L}}^{\otimes n})$, $n\geq 1$, forment
une base de la topologie de Zariski de $\widetilde{P}$.

\end{defin}

\vspace{-2mm}
\begin{defin}
Une {\it {\'e}toile} $\Sigma$ de $\mathfrak{a}$ 
est un sous-ensemble fini de $\mathfrak{a}$ tel que 
$0\in \Sigma$, $\Sigma=-\Sigma$ et $\Sigma$ contient une base de $\mathfrak{a}$.
\end{defin}

Soit l'anneau gradu{\'e} :
$\mathcal{R}=\sum_{k=0}^{\infty} 
 K[\mathfrak{X}^{\alpha};\alpha\in \mathfrak{a}]\cdot
\theta^k.$

\smallskip

\noindent On  d{\'e}finit une action du groupe $\mathfrak{b}$ sur $\mathcal{R}$ par : 
$\begin{cases} T_{\beta}^*(c)=c$, pour  $ c\in K, \\
T_{\beta}^*(\mathfrak{X}^{\alpha})=
\mathfrak{X}^{\alpha}(\beta)\mathfrak{X}^{\alpha}$, pour  
$\alpha \in \mathfrak{a}, \\ 
T_{\beta}^*(\theta)=
\mathfrak{X}^{\phi(\beta)}(\beta)\mathfrak{X}^{2\phi(\beta)}\theta.
\end{cases}$

\begin{defin}
Soit $\Sigma$ une {\'e}toile de $\mathfrak{a}$; on note $R_{\phi,\Sigma}$ le sous anneau de $\mathcal{R}$
engendr{\'e} sur $R$ par les {\'e}l{\'e}ments $T_{\beta}^*(\mathfrak{X}^{\alpha}\theta)$ pour
$\beta\in\mathfrak{b}$ et $\alpha\in \Sigma$, i.e. :

$$R_{\phi,\Sigma}=
R[\mathfrak{X}^{\phi(\beta)+\alpha}(\beta)\mathfrak{X}^{2\phi(\beta)+\alpha}\theta
]_{\beta\in\mathfrak{b},\alpha\in \Sigma}.$$
\end{defin}

D'apr{\`e}s le lemme \ref{lemme2} on peut supposer, quitte {\`a} remplacer
$\phi$ par $n\phi$, que $R_{\phi,\Sigma}\subset 
R[\mathfrak{X}^{\alpha}\theta]_{\alpha\in \mathfrak{a}}$.

On montre alors que  $\Proj(R_{\phi,\Sigma})$ est 
un mod{\`e}le relativement complet pour $\widetilde{G}$. Comme 
$R_{\phi,\Sigma}$ est un anneau gradu{\'e} engendr{\'e} par ses 
{\'e}l{\'e}ments de degr{\'e} $1$, $\Proj(R_{\phi,\Sigma})$ est muni d'un faisceau
tr{\`e}s 
ample inversible canonique, qui est le $\mathcal{O}(1)$.

On obtient ainsi le :

\vspace{-4mm}
\begin{theo} $\mathrm{(Mumford\enspace \cite{mum})}$
Soit $\widetilde{G}$ un tore d{\'e}ploy{\'e} sur $S$, 
$\mathfrak{b}\subset\widetilde{G}_\eta$ un 
groupe de p{\'e}riodes et $\phi:\mathfrak{b}\rightarrow \mathfrak{a}$ une polarisation. Alors,
pour toute {\'e}toile $\Sigma$ de $\mathfrak{a}$, 
quitte {\`a} remplacer $\phi$ par $n\phi$ $(n\in\Z$, $n\gg 0)$, 
$\widetilde{P}=\Proj(R_{\phi,\Sigma})$, muni de  son faisceau canonique $\mathcal{O}(1)$, 
est un mod{\`e}le relativement complet pour
 $\widetilde{G}$ sur $S$, par rapport {\`a} $(\mathfrak{b},2\phi)$.
\end{theo}

On remarque que $\widetilde{G}_\eta=\widetilde{P}_\eta$.

La construction du quotient proc{\`e}de en deux temps : Mumford forme d'abord le  quotient 
$\mathfrak{P}$ du compl{\'e}t{\'e} formel de $\widetilde{P}$ le long du bord, par $\mathfrak{b}$. 
Ce quotient est un sch{\'e}ma formel projectif et de type fini, donc s'alg{\'e}brise en 
un sch{\'e}ma projectif de type fini not{\'e} $P$.

Consid{\'e}rons l'ouvert 
$\bigcup_{\beta\in\mathfrak{b}}T_{\beta}(\widetilde{G})\subset\widetilde{P}$.
Soit $\widetilde{B}=\widetilde{P}-\bigcup_{\beta\in\mathfrak{b}}T_{\beta}(\widetilde{G})$
le sous-sch{\'e}ma r{\'e}duit, 
et  $\mathfrak{B}$ le quotient par $\mathfrak{b}$ de son compl{\'e}t{\'e}  formel. 
C'est la compl{\'e}tion formelle d'un sous-sch{\'e}ma 
r{\'e}duit $B\subset P$. Posons $G=P\bs B$.
Par construction les compl{\'e}tions $I$-adiques de $G$ 
et  $\widetilde{G}$ sont canoniquement isomorphes.

\vspace{-2mm}
\begin{theo} $\mathrm{(Mumford \enspace\cite{mum})}$ Le  sch{\'e}ma $G/S$ est 
semi-ab{\'e}lien, 
 $G_\eta$  est une vari{\'e}t{\'e} ab{\'e}lienne  et $G_0$ est un tore d{\'e}ploy{\'e} de rang $r$.
 Le  sch{\'e}ma $G/S$ ne d{\'e}pend que du tore $\widetilde{G}$ et 
 du groupe de p{\'e}riodes $\mathfrak{b}$, et  il est ind{\'e}pendant de la
fonction de polarisation $\phi$ et  du mod{\`e}le relativement complet $\widetilde{P}$.
La construction de $G/S$ est fonctorielle en $\widetilde{G}/S$ et en $\mathfrak{b}$.
\end{theo}

\vspace{-.4cm} 
\section{Construction de VAHB d{\'e}g{\'e}n{\'e}rantes.}\label{vahb-degenerantes}

On applique la construction de Mumford pour construire des 
vari{\'e}t{\'e}s ab{\'e}liennes de Hilbert-Blumenthal d{\'e}g{\'e}n{\'e}rantes.
Afin d'{\'e}viter des r{\'e}p{\'e}titions avec la partie 2 de \cite{dimtildg}, 
nous n'allons  donner la d{\'e}finition d'une vari{\'e}t{\'e} ab{\'e}lienne de Hilbert-Blumenthal
que dans le cas  o{\`u} le discriminant $\Delta_F$ du 
corps $F$ est inversible.

\begin{defin} Une vari{\'e}t{\'e} ab{\'e}lienne de Hilbert-Blumenthal (abr{\'e}g{\'e} VAHB) 
 sur un $\Z[\frac{1}{\Delta_F}]$-sch{\'e}ma $S$ est la donn{\'e}e d'un sch{\'e}ma   
ab{\'e}lien  $f:A\rightarrow S$ de dimension relative $d_F$ et 
d'une injection  $\iota:\mathfrak{o}\hookrightarrow \End(A/S)$ tels
que le faisceau $\underline{\omega}=f_*\Omega^1_{A/S}$ soit   
localement libre de rang $1$ sur $\mathfrak{o}\otimes\mathcal{O}_S$,
pour la topologie de Zariski.
\end{defin}

Pour tout id{\'e}al fractionnaire $\mathfrak{f}$ de $F$
on pose $\mathfrak{f}^*=\mathfrak{f}^{-1}\mathfrak{d}^{-1}$.
On a un accouplement parfait $\Tr_{F\!/\!\Q} :
\mathfrak{f}\times\mathfrak{f}^* \rightarrow\Z$.
 
\medskip
Soit $X$ un id{\'e}al fractionnaire de $F$, muni de sa positivit{\'e} 
$X_+=X\cap(F\otimes \R)_+$.

\smallskip
\noindent{\bf L'anneau de base} $\overline{S}_{\sigma}$.
 Soit $R=\Z[q^\xi;\xi\in X]$.

Soit $S=\Spec(R)=\mathbb{G}_m\otimes X^{*}$ le tore 
de groupe de caract{\`e}res  $X$.

Soit  $\Sigma$ un {\'e}ventail complet lisse   
de $X_{+}^{*}$ et soit $S\hookrightarrow
S_{\Sigma}$, l'immersion torique associ{\'e}e. On rappelle
qu'elle est obtenue en recollant,  pour $\sigma\in \Sigma$, les immersions toriques
affines $S\hookrightarrow
S_{\sigma}=\Spec(R_{\sigma})$, o{\`u} 
$R_{\sigma}=\Z[q^\xi;\xi\in X\cap\check{\sigma}]$.
Soit $S_{\sigma}^{\wedge}$ le compl{\'e}t{\'e} de $S_{\sigma}$
le long de $S_{\sigma}^{\infty}:=S_{\sigma}\backslash S$ et 
$S_{\Sigma}^{\wedge}$ le compl{\'e}t{\'e} de $S_{\Sigma}$
le long de $S_{\Sigma}^{\infty}:=S_{\Sigma}\backslash S$.

Pour {\'e}crire les choses plus explicitement, donnons nous une 
base   $\xi_1^{*}$,..,$\xi_r^{*}$  de $\sigma$
que l'on compl{\`e}te en une base  $\xi_1^{*}$,..,$\xi_d^{*}$ de  $X^{*}$. 
Soit $\xi_1$,..,$\xi_d$ la base duale de $X$ et posons $Z_i=q^{\xi_i}$.
Alors $R_{\sigma}=\Z[Z_1,..,Z_r,Z_{r+1}^{\pm}Z_d^{\pm}]$ et 
$S_{\sigma}^{\infty}$
est le  diviseur {\`a} croisements normaux de $S_{\sigma}$
d{\'e}fini par l'{\'e}quation $Z_1...Z_r=0$.

On a $S_{\sigma}^{\wedge}=\Spf(R_{\sigma}^{\wedge})$,
o{\`u} $R_{\sigma}^{\wedge}$ est le compl{\'e}t{\'e} de 
$R_{\sigma}$ en l'id{\'e}al principal radiciel $(Z_1\cdot ...\cdot Z_r)$.

Pour d{\'e}crire ce compl{\'e}t{\'e}, on d{\'e}compose tout 
$\underline{n}=(n_1,..,n_d)\in\Z^d$ en 
$(\underline{n}',\underline{n}'')\in \Z^r\times \Z^{d-r}$. 
Disons qu'une s{\'e}rie de Laurent formelle $\sum_{\underline{n}\in\Z^d} 
c_{\underline{n}} Z_1^{n_1}...Z_d^{n_d}$ {\`a} coefficients $c_{\underline{n}}\in\Z$
est $(Z_1\cdot\ldots \cdot Z_r)$-enti{\`e}re si

\noindent (i) pour tout $\underline{n}''$, 
$c_{\underline{n}',\underline{n}''}=0$, si $\underline{n}'\not\in \n^r$,

\noindent (ii)  pour tout $H\geq 1$ on a 
$c_{\underline{n}',\underline{n}''}=0$, 
pour presque tout $(\underline{n}',\underline{n}'')\notin {[H,\infty[}^r\times \Z^{d-r}$.

Le compl{\'e}t{\'e}
$ R_{\sigma}^{\wedge}$ s'identifie alors {\`a} l'ensemble des s{\'e}ries 
$\sum_{\underline{n}\in\Z^d} 
c_{\underline{n}} Z_1^{n_1}...Z_d^{n_d}$ qui sont 
$(Z_1\cdot\ldots\cdot Z_r)$-enti{\`e}res. C'est un anneau normal.

On voit ainsi que  $R_{\sigma}^{\wedge}$ est aussi le compl{\'e}t{\'e} 
de $R_{\sigma}$ par rapport {\`a} la topologie suivante :

\vspace{-4mm}
\begin{equation}
q^{\xi_i}\rightarrow 0\iff \Tr_{F\!/\!\Q}(\xi_i \xi^*)\rightarrow +\infty, \enspace \forall 
\xi^*\in  \sigma.
\end{equation}

L'anneau de base  sur lequel nous effectuons la construction de Mumford ici 
est $R_{\sigma}^{\wedge}$. Soit 
$\overline{S}_{\sigma}=\Spec(R_{\sigma}^{\wedge})$; 
posons $\overline{S}{}^{0}_{\sigma}=
S\underset{S_{\sigma}}{\times}\overline{S}_{\sigma}=\Spec(R_{\sigma}^{\wedge}
\otimes_{R_{\sigma}}R)$. C'est l'ouvert de
$\overline{S}_{\sigma}$ obtenu en rendant inversible  $q^{\xi}$ 
pour tout {\'e}l{\'e}ment $\xi$ de $X\cap \check\sigma^0$ 
(o{\`u} $\check\sigma^0$ d{\'e}signe l'int{\'e}rieur du c{\^o}ne dual 
$\check\sigma$ de $\sigma$). Soit $\overline{S}_{\sigma 0}:=
\overline{S}_{\sigma}\bs\overline{S}{}^0_{\sigma}$
muni de la structure r{\'e}duite. Si $\sigma'\subset\sigma$, on a une fl{\`e}che 
$\overline{S}_{\sigma'} \rightarrow \overline{S}_{\sigma}$.

\smallskip

\noindent{\bf Le tore $\widetilde{G}$.}
Soit $\mathfrak{a}$ (=$P^*$ dans les notations de Rapoport \cite{rapoport}) 
un id{\'e}al du corps de nombres totalement r{\'e}el $F$ et soit 
$\widetilde{G}:=(\mathbb{G}_m\otimes \mathfrak{a}^*)\times \overline{S}_{\sigma}$ 
le $\overline{S}_{\sigma}$-tore de groupe des caract{\`e}res $\mathfrak{a}$. 
Explicitement :
$\widetilde{G}=\Spec\left(
R_{\sigma}^{\wedge}[\mathfrak{X}^\alpha ;\alpha\in \mathfrak{a}]\right)$.

\smallskip
\noindent{\bf L'ensemble des p{\'e}riodes $\mathfrak{b}$.}
Soit $\mathfrak{b}$ (=$N$ dans les notations de Rapoport \cite{rapoport}) 
un id{\'e}al fractionnaire de $F$, tel que 

$$  \mathfrak{ab}^{-1}=\mathfrak{c} \text{ et }
\mathfrak{ab}\subset X$$

Pour chaque $\beta\in \mathfrak{b}$ on d{\'e}finit un 
$\overline{S}{}^{0}_{\sigma}$-point de 
$\widetilde{G}$, par le morphisme 

$$R_{\sigma}^{\wedge}[
\mathfrak{X}^\alpha;\alpha\in \mathfrak{a}]\rightarrow R_{\sigma}^{\wedge}
\otimes_{R_{\sigma}}R, \enspace
\mathfrak{X}^\alpha\mapsto q^{\alpha \beta}.$$

Ceci d{\'e}finit un homomorphisme $\mathfrak{o}$-{\'e}quivariant 
de $\overline{S}{}^{0}_{\sigma}$-sch{\'e}mas en groupes \\
$q : \mathfrak{b}\rightarrow \mathbb{G}_m\otimes \mathfrak{a}^*=\widetilde{G}$,
(o{\`u} $\mathfrak{b}$ d{\'e}signe le sch{\'e}ma en groupes constant).

\smallskip
\noindent{\bf La polarisation} $\phi$.
Se donner une polarisation $\mathfrak{o}$-lin{\'e}aire $\phi:\mathfrak{b}\rightarrow 
\mathfrak{a}$
(voir la d{\'e}finition \ref{polar}) revient {\`a} se donner un {\'e}l{\'e}ment 
$[\phi]\in \mathfrak{c}_{+}=\mathfrak{c}\cap(F\otimes \R)_+$.

\smallskip

La construction de Mumford donne un sch{\'e}ma semi-ab{\'e}lien 
$G_{\sigma}$ sur $\overline{S}_{\sigma}$.

\smallskip

\noindent{\bf Propri{\'e}t{\'e}s du sch{\'e}ma semi-ab{\'e}lien $G_{\sigma}$.} 

\noindent  $\relbar$ La restriction de $G_{\sigma}$ {\`a} $\overline{S}{}^{0}_{\sigma}$
est une VAHB, not{\'e}e $G_{\sigma}^0$.

\smallskip
\noindent  $\relbar$ Tout {\'e}l{\'e}ment $[\phi]\in \mathfrak{c}$ donne une fl{\`e}che 
naturelle $\mathbb{G}_m\otimes \mathfrak{a}^*\rightarrow 
\mathbb{G}_m\otimes \mathfrak{b}^*$, d'o{\`u}, par fonctorialit{\'e} de 
la construction, une fl{\`e}che sym{\'e}trique $\phi$ de la vari{\'e}t{\'e} ab{\'e}lienne 
 $G_{\sigma}^0= (\mathbb{G}_m\otimes \mathfrak{a}^*)/q(\mathfrak{b})$
 vers sa duale  $(G_{\sigma}^0)^t=
(\mathbb{G}_m\otimes \mathfrak{b}^*)/q(\mathfrak{a})$. 
Si $[\phi]\in \mathfrak{c}_{+}$, alors $\phi$  est une polarisation.

\smallskip
\noindent  $\relbar$ Par le lemme du serpent, appliqu{\'e} {\`a} la 
multiplication par $\mathfrak{n}$ 
dans $\Gm\otimes \mathfrak{a}^*$, on trouve la $\mathfrak{n}$-torsion de 
$G_{\sigma}^0$ (qui est le sous-sch{\'e}ma  en groupes r{\'e}duit,  
intersection des noyaux des multiplications par les {\'e}l{\'e}ments de $\mathfrak{n}$)
au milieu de la suite exacte 
\begin{equation}\label{ntorsion}
1 \rightarrow(\mathfrak{a}/\mathfrak{na})(1) \rightarrow 
G_{\sigma}^0[\mathfrak{n}]
\rightarrow \mathfrak{n}^{-1}\mathfrak{b}/\mathfrak{b} \rightarrow 0
\end{equation}

\smallskip
\noindent  $\relbar$ La restriction de $G_{\sigma}$ {\`a} 
$\overline{S}_{\sigma 0}$
est {\'e}gale au tore $(\Gm\otimes \mathfrak{a}^*)\times\overline{S}_{\sigma 0}$.

\smallskip
\noindent $\relbar$ La construction est  fonctorielle en les $\sigma\in\Sigma$
et compatible avec l'action de $\mathfrak{o}^{\times}$,
i.e.   pour  tout  $\sigma'\subset\sigma$ et pour tout
$u\in \mathfrak{o}^{\times}$ on a des   diagrammes  cart{\'e}siens :

\vspace{.2cm} 

$ $\hspace{2cm}
\xymatrix@R=10pt{G_{\sigma'}\ar[d]\ar[r] & G_{\sigma}\ar[d]\\
\overline{S}_{\sigma'}\ar[r]  & \overline{S}_{\sigma}}\hspace{2cm}
\xymatrix@R=10pt{G_{\sigma}\ar[d]\ar^{\sim}[r] & G_{u^2\sigma}\ar[d]\\
\overline{S}_{\sigma}\ar^{\sim}[r]  & \overline{S}_{u^2\sigma}}

\section{$R$-pointes et $(R,\mathfrak{n})$-pointes.}

Pour tout id{\'e}al $\mathfrak{f}\subset\mathfrak{o} $ on note 
$\mathfrak{o}_{\mathfrak{f}}^\times$ le sous-groupe de 
$\mathfrak{o}^\times$ form{\'e} des unit{\'e}s congrues {\`a} $1$ modulo $\mathfrak{f}$. 
On note $\mathfrak{o}_+^\times$ le groupe des unit{\'e}s totalement 
positives de $\mathfrak{o}$.

 Pour tout $\mathfrak{o}$-r{\'e}seau $L$ de $F^2$ notons 
$G^+(L)$ le stabilisateur de $L$ dans $G^D_{\Q+}$ 
(pour l'action {\`a} gauche donn{\'e}e par 
$\gamma\cdot l=l\gamma^{-1}$, pour tout 
 $\gamma\in G^D_{\Q}$ et $l\in L$). On a 

$$G^+(\mathfrak{o}\oplus\mathfrak{c}^*)=
\Big{\{}\gamma \in \begin{pmatrix}\mathfrak{o} &\mathfrak{c}^* \\
\mathfrak{cd} &\mathfrak{o} \end{pmatrix}\Big{|} 
\nu(\gamma)\in \mathfrak{o}_+^\times\cap D_{\Q}\Big{\}}.$$ 

$$\text{Posons }\Gamma=\Gamma^D_1(\mathfrak{c},\mathfrak{n})
=\Big{\{} \begin{pmatrix}a &b \\c &d\end{pmatrix}
\in G^+(\mathfrak{o}\oplus\mathfrak{c}^*)
\enspace \Big{|} \enspace c\in \mathfrak{cdn},  
\enspace  d\equiv 1 (\Mod \mathfrak{n}) \Big{\}}.$$

Cette partie {\'e}tudie la combinatoire des pointes d'une vari{\'e}t{\'e} 
modulaire de Hilbert-Blumenthal en niveau 
$\Gamma_1^D(\mathfrak{c},\mathfrak{n})$ 
et servira {\`a}  la construction 
de cartes locales pour les compactifications toro{\"\i}dales.
Cette {\'e}tude a {\'e}t{\'e} d{\'e}j{\`a} effectu{\'e}e par Rapoport en 
niveau $\Gamma^D(\mathfrak{c},\mathfrak{n})$ et en niveau 
$\Gamma_1^D(\mathfrak{c},\mathfrak{n})$ pour une pointe 
non-ramifi{\'e}e,  lorsque $\mathfrak{n}$ est un entier naturel et $D=\Gm$  
(voir  \cite{rapoport}). 
Par ailleurs, lorsque $F=\Q$,  l'{\'e}tude est faite par
Deligne et Rapoport \cite{DeRa}, en niveau $\Gamma(\mathfrak{n})$,
et  par Katz et Mazur \cite{KaMa} en g{\'e}n{\'e}ral.

Soit $\mathfrak{c}$  un id{\'e}al fractionnaire de $F$, muni de sa
 de positivit{\'e}  $\mathfrak{c}_+=\mathfrak{c}\cap(F\otimes\R)_+$.

Les objets combinatoires consid{\'e}r{\'e}s dans cette partie sont inspir{\'e}s
par les  structures de niveau des VAHB :
une VAHB $\mathfrak{c}$-polaris{\'e}e complexe admet une uniformisation 
de la forme $F\otimes \C/L$, o{\`u} $L$ est un $\mathfrak{o}$-r{\'e}seau de
$F^2$ tel que  $\wedge^2_{\mathfrak{o}}L=\mathfrak{c}^*$. Or, un 
tel r{\'e}seau s'{\'e}crit $L=\mathfrak{b}\oplus\mathfrak{a}^*$, avec 
$\mathfrak{a}$ et $\mathfrak{b}$ deux id{\'e}aux fractionnaires de $F$
tels que $\mathfrak{a}^*\mathfrak{b}= \mathfrak{c}^*$.
La $\mu_{\mathfrak{n}}$-structure de niveau sur une telle VAHB est donn{\'e}e 
alors par un homomorphisme injectif de  $\mathfrak{o}$-modules  
$\beta : \mathfrak{n}^{-1}\mathfrak{d}^{-1}/\mathfrak{d}^{-1} 
\hookrightarrow \mathfrak{n}^{-1}L/L$. Par ailleurs
 tout $\mathfrak{o}$-module projectif 
de rang $2$ est isomorphe {\`a} un $\mathfrak{o}$-r{\'e}seau de $F^2$. 
La d{\'e}finition suivante est une variante de celle donn{\'e}e par Rapoport
dans le cas $D=\Gm$ :

\begin{defin} \label{R-pointe} Une  $R$-pointe $\mathcal{C}$ (resp. une 
classe d'isomorphisme de $R$-pointes) est une classe  d'{\'e}quivalence de 
 sextuplets  $(\mathfrak{a},\mathfrak{b},L,i,j,\lambda)$, o{\`u}

$\mathrm{(i)}$ $\mathfrak{a}$ et $\mathfrak{b}$ sont deux id{\'e}aux 
fractionnaires de $F$ tels que $\mathfrak{a}^*\mathfrak{b}=\mathfrak{c}^*$,

$\mathrm{(ii)}$ $L$ est un $\mathfrak{o}$-r{\'e}seau de $F^2$ tel que
 l'on a une suite  exacte $\mathfrak{o}$-modules 
$$ 0\rightarrow \mathfrak{a}^* \stackrel{i}
{\rightarrow} L \stackrel{j}{\rightarrow} \mathfrak{b} \rightarrow 0,$$

$\mathrm{(iii)}$  $\lambda: \wedge^2_{\mathfrak{o}}L \rightarrow \mathfrak{c}^*$
est un isomorphisme $\mathfrak{o}$-lin{\'e}aire (polarisation),

\noindent pour la relation d'{\'e}quivalence suivante :
$(\mathfrak{a},\mathfrak{b},L,i,j,\lambda)$  et 
$(\mathfrak{a'},\mathfrak{b'},L',i',j',\lambda')$ sont {\'e}quivalents, 
si $\mathfrak{a}=\mathfrak{a'}$,
 $\mathfrak{b}=\mathfrak{b'}$ (resp. $\mathfrak{a}=\xi\mathfrak{a'}$
et $\mathfrak{b}=\xi\mathfrak{b'}$ avec $\xi\in F$) 
et s'il existe un diagramme commutatif 
de $\mathfrak{o}$-modules :

\hspace{2cm}\xymatrix@R=15pt{ 0 \ar[r] & \mathfrak{a}^* \ar^{i}[r]\ar[d] &
L \ar^{j}[r]\ar[d] & \mathfrak{b}\ar[r]\ar[d] & 0,\\
 0 \ar[r] & \mathfrak{a'}^* \ar^{i'}[r] &
L' \ar^{j'}[r] & \mathfrak{b'}\ar[r] & 0}

\noindent o{\`u} les fl{\`e}ches verticales sont des isomorphismes et tel que 
l'isomorphisme $\wedge^2_{\mathfrak{o}}L \cong
\wedge^2_{\mathfrak{o}}L'$ (d{\'e}duit de $L \cong L'$) 
induise, via $\lambda$ et  $\lambda'$, un automorphisme de $\mathfrak{c}^*$,
donn{\'e} par un {\'e}l{\'e}ment de $\mathfrak{o}_{D+}^\times:=\mathfrak{o}_+^\times\cap D_{\Q}$.
\end{defin}

L'application qui {\`a} une  $R$-pointe 
$\mathcal{C}=(\mathfrak{a},\mathfrak{b},L,i,j,\lambda)$ associe 
l'id{\'e}al $\mathfrak{b}$ est une bijection entre l'ensemble 
des $R$-pointes et l'ensemble $\mathcal{I}_F$ des id{\'e}aux fractionnaires
de $F$. En effet, par (i) la donn{\'e}e de $\mathfrak{b}$
d{\'e}termine $\mathfrak{a}=\mathfrak{bc}$, et 
deux suites exactes courtes (ii),  correspondant au m{\^e}me
 id{\'e}al $\mathfrak{b}$, sont {\'e}quivalentes, car toutes les deux sont 
scind{\'e}es. 

La notion d'isomorphisme de $R$-pointes  correspond alors {\`a}  celle 
d'homoth{\'e}tie des id{\'e}aux. On obtient par  passage  au quotient un 
isomorphisme entre les classes d'isomorphisme
de $R$-pointes et le groupe  $\mathrm{Cl}_F$ des classes d'id{\'e}aux de $F$.

Une $R$-pointe est d{\'e}termin{\'e}e par son $\mathfrak{o}$-r{\'e}seau $L$ de $F^2$
 (en effet, la donn{\'e}e d'un tel r{\'e}seau d{\'e}termine  les id{\'e}aux 
$\mathfrak{a}^*:=L\cap(\{0\} \times F)$ et $\mathfrak{b}=\mathfrak{ca}^{-1}$, 
et donc la $R$-pointe $\mathcal{C}$, {\`a} {\'e}quivalence pr{\`e}s). 
Le groupe $ G_{\Q}^{\mathfrak{o}}:=\{\gamma\in G^D_{\Q}| 
\nu(\gamma)\in \mathfrak{o}_+^\times\}$ 
agit transitivement sur ces r{\'e}seaux. 
Le stabilisateur du r{\'e}seau $\mathfrak{o}\oplus \mathfrak{c^*}$ dans $ G_{\Q}^{\mathfrak{o}}$
est {\'e}gal {\`a}  $G^+(\mathfrak{o}\oplus\mathfrak{c}^*)$. De plus, deux
r{\'e}seaux $L$ et $L'$ donnent la m{\^e}me $R$-pointe $\mathcal{C}$, 
si et seulement s'ils sont dans la m{\^e}me 
$T_{\Z}U_{\Q}$-orbite. 
Le  diagramme commutatif suivant, traduit la correspondance entre 
les $R$-pointes et les pointes classiques dans $\p^{1}(F)$ pour
le sous-groupe de congruence $G^+(\mathfrak{o}\oplus\mathfrak{c}^*)$

$$\xymatrix@R=10pt@C=15pt{ 
 \mathcal{I}_F  \ar^{{}\hspace{-5mm}\sim}[r]\ar[d] &  
R\mathrm{-pointes} \ar^{{}\hspace{-10mm}\sim}[r]\ar[d] & 
T_{\Z}U_{\Q} \bs G_{\Q}^{\mathfrak{o}}/G^+(\mathfrak{o}\oplus\mathfrak{c}^*)
\ar^{{}\hspace{-2mm}\sim}[r]\ar[d]  & 
G^+(\mathfrak{o}\oplus\mathfrak{c}^*)\bs
F^2-\{0\}/\mathfrak{o}^\times \ar[d]\\
\mathrm{Cl}_F \ar^{{}\hspace{-10mm}\sim}[r] &  
R\mathrm{-pointes/isom.}   \ar^{{}\hspace{-5mm}\sim}[r] & 
 B_{\Q}\bs  G^D_{\Q}/ G^+(\mathfrak{o}\oplus\mathfrak{c}^*) \ar^{{}\hspace{2mm}\sim}[r]&
G^+(\mathfrak{o}\oplus\mathfrak{c}^*)\bs \p^1(F),}$$

\noindent o{\`u} pour tout $\gamma= \begin{pmatrix} a & b \\ c& d\end{pmatrix} 
\in G_{\Q}^{\mathfrak{o}}$  la double classe 
$B_{\Q}\gamma^{-1} G^+(\mathfrak{o}\oplus\mathfrak{c}^*)$ s'envoie 
d'une part sur la pointe classique 
$G^+(\mathfrak{o}\oplus\mathfrak{c}^*) \gamma \infty$ et d'autre part sur 
l'id{\'e}al $\mathfrak{b}=a\mathfrak{o}+c\mathfrak{c}^*$ (voir \cite{dimtildg} Lemme 1.7). 

\vspace{-4mm}
\begin{defin} \label{Rn-pointe}

$\mathrm{(i)}$ Une  $(R,\mathfrak{n})$-pointe $\mathcal{C}$ (resp. une 
classe d'isomorphisme de $(R,\mathfrak{n})$-pointes) est la donn{\'e}e
d'une classe d'{\'e}quivalence de paires form{\'e}es 
d'un sextuplet $(\mathfrak{a},\mathfrak{b},L,i,j,\lambda)$ 
(comme dans la d{\'e}finition \ref{R-pointe}) 
et d'un morphisme injectif  de  $\mathfrak{o}$-modules 

\vspace{-2mm}
 $$\beta : \mathfrak{n}^{-1}\mathfrak{d}^{-1}/\mathfrak{d}^{-1}
 \hookrightarrow \mathfrak{n}^{-1}L/L,$$
pour  la relation d'{\'e}quivalence suivante:

\smallskip
\noindent $\mathcal{C}$ est {\'e}quivalent {\`a}  $\mathcal{C}'$, s'il existe un 
isomorphisme de $\mathfrak{o}$-modules  
$L\cong L'$ induisant une {\'e}galit{\'e} (resp. un isomorphisme) 
des $R$-pointes sous-jacentes et  dont la r{\'e}duction 
modulo $\mathfrak{n}$ rend le diagramme suivant commutatif   :

\vspace{-2mm}
$$\xymatrix@R=12pt{ 
\mathfrak{n}^{-1}L/L\ar^{\sim}[rr] & & \mathfrak{n}^{-1}L'/L'. \\ 
& \mathfrak{n}^{-1}\mathfrak{d}^{-1}/\mathfrak{d}^{-1} \ar@{_{(}->}^{\beta}[ul] 
\ar@{^{(}->}_{\beta'}[ru] &}$$ 

On  associe {\`a}  $\mathcal{C}$  l'id{\'e}al  fractionnaire  
$\mathfrak{b}'\supset \mathfrak{b}$ tel que  
$\mathfrak{b}'/\mathfrak{b}= j(\im(\beta))$.

$\mathrm{(ii)}$
 Une $(R,\mathfrak{n})$-pointe  est dite {\it non-ramifi{\'e}e} lorsque
la fl{\`e}che $\beta:\mathfrak{n}^{-1}\mathfrak{d}^{-1}/\mathfrak{d}^{-1}
\hookrightarrow \mathfrak{n}^{-1}L/L$ se  factorise par la fl{\`e}che naturelle 
$\mathfrak{n}^{-1}\mathfrak{a}^*/\mathfrak{a}^* 
\hookrightarrow \mathfrak{n}^{-1}L/L$ (ou si de mani{\`e}re {\'e}quivalente 
$\mathfrak{b}'=\mathfrak{b}$).

$\mathrm{(iii)}$
 Soit une $(R,\mathfrak{n})$-pointe $\mathcal{C}$ et soit
 $n$ l'exposant du groupe  $\mathfrak{b'}/\mathfrak{b}$.
Une  $(R,\mathfrak{n})$-pointe $\mathcal{C}'$ est dite appartenir {\`a} la m{\^e}me 
{\it $(R,\mathfrak{n})$-composante} que $\mathcal{C}$ (resp. {\`a} une 
$(R,\mathfrak{n})$-composante isomorphe), 
s'il existe $\overline{a}\in (\Z/n)^\times$
et  un  isomorphisme de $\mathfrak{o}$-modules $L\cong L'$ induisant une {\'e}galit{\'e} 
(resp. un isomorphisme) des $R$-pointes sous-jacentes et  dont la r{\'e}duction $\psi$
modulo $\mathfrak{n}$ fait commuter le diagramme suivant 

$$\xymatrix@R=12pt{ \mathfrak{n}^{-1}L/L
\ar^{\sim}_{\varphi}[r] & 
\mathfrak{n}^{-1}L/L\ar^{\sim}_{\psi}[r] &  \mathfrak{n}^{-1}L'/L', \\ 
& \mathfrak{n}^{-1}\mathfrak{d}^{-1}/\mathfrak{d}^{-1} \ar@<3pt>@{_{(}->}^{\beta}[lu] 
\ar@<-2pt>@{^{(}->}_{\beta'}[ru] &  } $$

\noindent o{\`u} la fl{\`e}che $\varphi$ est un automorphisme $\mathfrak{o}$-lin{\'e}aire de 
$\mathfrak{n}^{-1}L/L$, induisant l'identit{\'e} sur 
$\mathfrak{n}^{-1}\mathfrak{a}^*/\mathfrak{a}^*$ et la multiplication par
$\overline{a}$ sur $\mathfrak{n}^{-1}\mathfrak{b}/\mathfrak{b}$.
\end{defin}

Soit $y_0$ tel que $\mathfrak{o}=\mathfrak{n}+y_0\mathfrak{c}$. 
On munit la $R$-pointe $L_0=\mathfrak{o}\oplus \mathfrak{c^*}$
de la structure de niveau $\beta_0:\mathfrak{n}^{-1}\mathfrak{d}^{-1}/\mathfrak{d}^{-1}
\overset{\cdot y_0}{\longrightarrow}
\mathfrak{n}^{-1}\mathfrak{c}^*/\mathfrak{c}^*\hookrightarrow \mathfrak{n}^{-1}L_0/L_0$.
Le groupe $ G_{\Q}^{\mathfrak{o}}$ agit transitivement sur ces r{\'e}seaux 
munis de structures de niveau et 
le stabilisateur de $(L_0,\beta_0)$ 
est $\Gamma$. De  plus, deux 
r{\'e}seaux $L$ et $L'$ donnent la m{\^e}me $R$-pointe $\mathcal{C}$, 
si et seulement s'ils sont dans la m{\^e}me 
$T_{\Z}U_{\Q}$-orbite. D'o{\`u} le diagramme suivant :

$$\xymatrix@R=12pt{ (R,\mathfrak{n})\mathrm{-pointes} \ar^{{}\hspace{-1mm}\sim}[r]\ar[d] &   
T_{\Z}U_{\Q}\bs G_{\Q}^{\mathfrak{o}}/\Gamma.\ar[d]\\
(R,\mathfrak{n})\mathrm{-pointes/isom.} \ar^{{}\hspace{5mm}\sim}[r] &  
B_{\Q}\bs G^D_{\Q}/\Gamma }$$

\begin{prop}\label{deuxideaux} $ $
 
 Soit une $(R,\mathfrak{n})$-pointe $\mathcal{C}$, donn{\'e}e par 
  $T_{\Z}U_{\Q}\gamma^{-1}\Gamma$,
$\gamma =\begin{pmatrix} a& b\\ c &d \end{pmatrix}\in G_{\Q}^{\mathfrak{o}}$.
Alors,

\noindent 
$\mathrm{(i)}$  L'id{\'e}al $\mathfrak{b}$, correspondant {\`a} la $R$-pointe 
sous-jacente {\`a} $\mathcal{C}$ est donn{\'e} par $a\mathfrak{o}+c\mathfrak{c}^* $ et 
sa classe  ne d{\'e}pend que de la classe d'isomorphisme de la pointe $\mathcal{C}$.

\smallskip
Quitte {\`a} changer $\gamma$, en le multipliant par un {\'e}l{\'e}ment de $U_{\Q}$, ce qui 
ne change pas sa classe double, on suppose que $\gamma\in G_{\Q}^{\mathfrak{o}}
\cap \begin{pmatrix} 
\mathfrak{b} & (\mathfrak{bc})^*\\ \mathfrak{bcd} & \mathfrak{b}^{-1}\end{pmatrix}$.
Sous cette hypoth{\`e}se:

\smallskip
\noindent 
$\mathrm{(ii)}$ La structure de niveau de  $\mathcal{C}$ est donn{\'e}e par 
$\beta : \mathfrak{n}^{-1}\mathfrak{d}^{-1}/\mathfrak{d}^{-1}
\overset{(y_0c,y_0d)}{\longrightarrow} \mathfrak{n}^{-1}L/L$, o{\`u} $L=\mathfrak{b}\oplus
\mathfrak{a}^*$ , avec $\mathfrak{a}=\mathfrak{bc}$.

\smallskip
\noindent 
$\mathrm{(iii)}$ L'id{\'e}al 
 $\mathfrak{b}'$ de la d{\'e}finition \ref{Rn-pointe}(i) est contenu dans  $\mathfrak{n}^{-1} \mathfrak{b}$ et 
sa classe  ne d{\'e}pend que de la classe d'isomorphisme de la pointe $\mathcal{C}$.
De plus  $\mathfrak{b'}=a\mathfrak{o}+c(\mathfrak{cn})^* $.
 La pointe $\mathcal{C}$ est non-ramifi{\'e}e, si et seulement si,
$c\in \mathfrak{nbcd}$.

\smallskip
\noindent 
$\mathrm{(iv)}$ Le groupe d'automorphismes de la $(R,\mathfrak{n})$-pointe $\mathcal{C}$
est {\'e}gal {\`a} $\gamma^{-1}\Gamma\gamma\cap B_{\Q}$.  La suite exacte 
$1 \rightarrow U \rightarrow B \rightarrow T \rightarrow 1$, donne une suite exacte :

$$ 0 \rightarrow X^*  \rightarrow \gamma^{-1}\Gamma\gamma\cap B_{\Q}  \rightarrow 
 \mathfrak{o}_{\mathcal{C}}^\times \rightarrow 1,$$
o{\`u} $X=\mathfrak{cbb'}$ et 
$\mathfrak{o}_{\mathcal{C}}^\times=
\{(u,\epsilon)\in \mathfrak{o}^\times \times \mathfrak{o}_{D+}^\times \quad|\quad
u-1\in \mathfrak{nb'b}^{-1},\quad u\epsilon-1\in \mathfrak{bb'}^{-1}\}$.
En particulier, on a  
$\mathfrak{o}_{\mathcal{C},1}^\times:=\mathfrak{o}_{\mathcal{C}}^\times\cap T_1=
\{u\in \mathfrak{o}^\times | u\in (1+\mathfrak{bb'}^{-1})
\cap (1+\mathfrak{nb'b}^{-1})\}.$

\smallskip
\noindent 
$\mathrm{(v)}$ L'ensemble des  $(R,\mathfrak{n})$-pointes est fibr{\'e} au-dessus de 
 $\mathcal{I}_F$.
La fibre de l'id{\'e}al $\mathfrak{b}$  est isomorphe {\`a}
$ (G^+(\mathfrak{b}\oplus\mathfrak{a}^*)\cap T_{\Z}U_{\Q})\bs
G^+(\mathfrak{b}\oplus\mathfrak{a}^*)
/\gamma^{-1}\Gamma\gamma $, o{\`u} 
$\mathfrak{a}=\mathfrak{bc}$,  $L=\mathfrak{b}\oplus\mathfrak{a}^*$.
Elle s'identifie avec l'ensemble :
$$(\mathfrak{n}^{-1}L/L)_{\mathrm{prim}}\Big/\Big\{
\begin{pmatrix} u \epsilon& \xi^*\\ 0 & u^{-1}\end{pmatrix}
\Big{|} u\in \mathfrak{o}^\times,\epsilon\in \mathfrak{o}_{D+}^\times, 
\xi^* \in (\mathfrak{cb^2})^*\Big\}, $$

\noindent 
o{\`u} $(\mathfrak{n}^{-1}L/L)_{\mathrm{prim}}$ d{\'e}signe l'ensemble des vecteurs 
primitifs du  $\mathfrak{o}/\mathfrak{n}$-module $\mathfrak{n}^{-1}L/L$,
et son  cardinal est {\'e}gal {\`a} 
$\underset{\mathfrak{n}^{-1}\mathfrak{b}\supset \mathfrak{b'}
\supset \mathfrak{b}}{\sum}
\#(\mathfrak{o}/\mathfrak{bb'}^{-1})^\times
\#(\mathfrak{o}/\mathfrak{nb'b}^{-1})^\times/
[(\mathfrak{o}^\times\times \mathfrak{o}_{D+}^\times):\mathfrak{o}_{\mathcal{C}}^\times].$

\smallskip
\noindent 
$\mathrm{(vi)}$ L'ensemble des  $(R,\mathfrak{n})$-composantes est 
fibr{\'e} au-dessus de   $\mathcal{I}_F$. La fibre de l'id{\'e}al $\mathfrak{b}$  
 s'identifie avec   l'ensemble :
$$(\mathfrak{n}^{-1}L/L)_{\mathrm{prim}}\Big/\Big\{
\begin{pmatrix} \overline{a} u \epsilon& \xi^* \\ 0 & u^{-1}\end{pmatrix}
\Big{|} u\in \mathfrak{o}^\times,  \epsilon\in \mathfrak{o}_{D+}^\times, 
\overline{a} \in (\Z/n)^\times
, \xi^*\in (\mathfrak{cb^2})^*\Big\}$$

\noindent 
qui est de   cardinal 
$\underset{\mathfrak{n}^{-1}\mathfrak{b}\supset \mathfrak{b'}\supset 
\mathfrak{b}}{\sum}\!\! \#(\mathfrak{o}/\mathfrak{bb'}^{-1})^\times
\#(\mathfrak{o}/\mathfrak{nb'b}^{-1})^\times/ \#(\Z/n)^\times
[(\mathfrak{o}^\times\times \mathfrak{o}_{D+}^\times)\!:\!
\mathfrak{o}_{\overline{\mathcal{C}}}^\times]$,  o{\`u} $n$ est  {\'e}gal {\`a}  
l'exposant   du groupe $\mathfrak{b'}/\mathfrak{b}$. De plus
 $$\mathfrak{o}_{\overline{\mathcal{C}}}^\times=
 \{(u,\epsilon)\in \mathfrak{o}^\times \times \mathfrak{o}_{D+}^\times 
\quad|\quad u-1\in \mathfrak{nb'b}^{-1},\quad u\epsilon\in 
(\Z/n)^\times+\mathfrak{bb'}^{-1}\},$$
$$\mathfrak{o}_{\overline{\mathcal{C}},1}^\times=
 \{ u \in \mathfrak{o}^\times 
\quad|\quad u\in (1+\mathfrak{nb'b}^{-1})\cap 
((\Z/n)^\times+\mathfrak{bb'}^{-1})\}.$$
\end{prop}

\noindent {\bf D{\'e}monstration : }
 (i) La $R$-pointe sous-jacente {\`a} $\mathcal{C}$ correspond {\`a} la 
 classe double $T_{\Z}U_{\Q}\gamma^{-1}G^+(\mathfrak{o}\oplus\mathfrak{c}^*)$
 et donc {\`a} la $G^+(\mathfrak{o}\oplus\mathfrak{c}^*)$-pointe 
 $\gamma\infty= \Big[\begin{matrix}a\\c\end{matrix}\Big]$. Par le diagramme
qui pr{\'e}c{\`e}de la d{\'e}finition \ref{Rn-pointe} la $R$-pointe 
 $\mathcal{C}$ correspond {\`a}  l'id{\'e}al   $\mathfrak{b}=a\mathfrak{o}+c\mathfrak{c}^* $.
  
 (ii)(iii) La structure de niveau $\beta$ de $L$ est obtenue en faisant agir
 $\gamma^{-1}$ sur la structure de niveau $\beta_0$ de  $L_0$.
 Or, par le choix que nous avons fait de $\gamma$, on a 
 $L_0\gamma=\mathfrak{b}\oplus \mathfrak{a}^*=L$ et donc 
 $\beta: \mathfrak{n}^{-1}\mathfrak{d}^{-1}/\mathfrak{d}^{-1}
\overset{ (cy_0,dy_0) }{\longrightarrow}\mathfrak{b}'/\mathfrak{b}\oplus
\mathfrak{n}^{-1}\mathfrak{a}^*/\mathfrak{a}^*\hookrightarrow \mathfrak{n}^{-1}L/L$.
 La pointe est donc non-ramifi{\'e}e si, et seulement, si  $cy_0\mathfrak{n}^{-1}\mathfrak{d}^{-1}
 \subset \mathfrak{b}$, i.e. $c\in \mathfrak{nbcd}$. Enfin 
$\mathfrak{b}'=\mathfrak{b}+cy_0 \mathfrak{d}^{-1}\mathfrak{n}^{-1}=
a\mathfrak{o}+c\mathfrak{c}^*+c\mathfrak{c}^*\mathfrak{n}^{-1}=
a\mathfrak{o}+c(\mathfrak{cn})^*$. L'ind{\'e}pendance des classes de
$\mathfrak{b}$ et $\mathfrak{b}'$ d{\'e}coule du lemme 1.7 de \cite{dimtildg} .

 (iv) Pour le calcul du groupe d'automorphismes $\gamma^{-1}\Gamma\gamma\cap B_{\Q}$ 
 de la $(R,\mathfrak{n})$-pointe $\mathcal{C}$, on remarque qu'il est form{\'e} de matrices
 $\begin{pmatrix} u\epsilon  & \xi_{u,\epsilon}^*\\  0 & u^{-1}\end{pmatrix}$, avec 
 $u\in\mathfrak{o}^\times$, $\epsilon \in\mathfrak{o}_{D+}^\times$, 
 $\xi^* \in (\mathfrak{cb}^2)^*$ (c'est la forme g{\'e}n{\'e}rale 
 d'un automorphisme de la   $R$-pointe sous-jacente) 
qui respectent en plus la structure de niveau $\beta$. Ceci {\'e}quivaut au syst{\`e}me
 \begin{equation}\label{Rn-niveau}
 \begin{cases} (u\epsilon-1)c\in  \mathfrak{nbcd}\\ 
 (u-1)d- \epsilon^{-1}\xi_{u,\epsilon}^*c\in 
 \mathfrak{ncda}^*=\mathfrak{nb}^{-1}\end{cases}.
\end{equation}

En posant $u=\epsilon=1$ on retrouve que $X^*$ est form{\'e} des
$\xi^*\in c^{-1}\mathfrak{nb}^{-1}\cap (\mathfrak{cb}^2)^*=
(\mathfrak{cb})^*((c(\mathfrak{cn}^*)^{-1}\cap \mathfrak{b}^{-1})=
(\mathfrak{cbb'})^*$, i.e. $X=\mathfrak{cbb'}$.

Pour le calcul de $\mathfrak{o}_{\mathcal{C}}^\times$ on remarque que
la premi{\`e}re condition de (\ref{Rn-niveau}) {\'e}quivaut {\`a} 
$u \epsilon-1\in c^{-1}\mathfrak{nbcd}\cap 
\mathfrak{o}=\mathfrak{b}(c(\mathfrak{cn})^{*-1}\cap \mathfrak{b}^{-1})=
\mathfrak{bb' }^{-1}$. 
La deuxi{\`e}me condition {\'e}quivaut {\`a} 
$u-1\in d^{-1}(\mathfrak{nb}^{-1}+c(\mathfrak{cb^2})^{*})=
(d\mathfrak{b})^{-1}\mathfrak{nb'b}^{-1}$. Par ailleurs 
$u-1\in \mathfrak{o}\subset c^{-1}\mathfrak{nb'cd}=(c(\mathfrak{bc})^*)^{-1}
\mathfrak{nb'b}^{-1}$. Comme 
$(d\mathfrak{b})^{-1}\cap (c(\mathfrak{bc})^*)^{-1}=
(d\mathfrak{b}+c(\mathfrak{bc})^*)^{-1}=\mathfrak{o}$, par le choix
de $\gamma$, on en d{\'e}duit que la deuxi{\`e}me condition de 
(\ref{Rn-niveau}) {\'e}quivaut {\`a} $u-1\in \mathfrak{nb'b}^{-1}$.

Notons  que pour tout $u\in \mathfrak{o}_{\mathcal{C}}^\times$, 
$\xi_{u,\epsilon}^*\in c^{-1}(\mathfrak{nb}^{-1}+d(\mathfrak{bb'}^{-1}\cap \mathfrak{nb'b}^{-1}))\subset (\mathfrak{cbb'})^*+d\mathfrak{b}((\mathfrak{cb}^2)^*\cap (\mathfrak{ncb'}^2)^*)
\subset (\mathfrak{cbb'})^*+(\mathfrak{cb}^2)^*\cap (\mathfrak{ncb'}^2)^*$, 
et ce dernier est un id{\'e}al inclus (parfois strictement!) 
dans $(\mathfrak{cb}^2)^*$ (voir l'exemple {\`a} la fin de l'article).

(v)(vi) Comme  $\gamma$ transforme
$\mathfrak{o}\oplus \mathfrak{c}^*$ en 
$\mathfrak{b}\oplus \mathfrak{a}^*$ et 
$\gamma^{-1} G^+(\mathfrak{o}\oplus \mathfrak{c}^*)\gamma=
G^+(\mathfrak{b}\oplus \mathfrak{a}^*)$, la fibre de l'id{\'e}al $\mathfrak{b}$  
est isomorphe {\`a} $ (G^+(\mathfrak{b}\oplus\mathfrak{a}^*)\cap T_{\Z}U_{\Q})\bs
G^+(\mathfrak{b}\oplus\mathfrak{a}^*)/\gamma^{-1}\Gamma\gamma $,
L'ensemble  $G^+(\mathfrak{b}\oplus\mathfrak{a}^*)/\gamma^{-1}\Gamma\gamma $
s'identifie avec celui des vecteurs 
primitifs du  $\mathfrak{o}/\mathfrak{n}$-module $\mathfrak{n}^{-1}L/L$.
Le calcul du cardinal de la fibre se fait en analysant la condition 
sous laquelle deux vecteurs primitifs correspondent {\`a} la m{\^e}me 
$(R,n)$-pointe. La d{\'e}monstration du (vi) est tout a fait analogue. 

Comme par d{\'e}finition $n\mathfrak{o}\subset \mathfrak{bb'}^{-1}
\subset \mathfrak{o}$, l'ensemble  $(\Z/n)^\times+\mathfrak{bb'}^{-1}$
est bien une r{\'e}union de classes de $\mathfrak{o}$, modulo l'id{\'e}al 
entier $\mathfrak{bb'}^{-1}$
Notons que 
$[\mathfrak{o}_{\overline{\mathcal{C}}}^\times:\mathfrak{o}_{\mathcal{C}}^\times]$
divise  $\#(\Z/n)^\times$  et le quotient repr{\'e}sente le nombre de 
$(R,n)$-pointes dans la $(R,n)$-composante $\overline{\mathcal{C}}$.
\hfill$\square$

\vspace{-3mm}
\begin{ex}\label{Rn-exemple} 
On pose  $\mathfrak{c}=\mathfrak{o}$ (polarisation principale)
et $G=G^*$ ($\mathfrak{o}_{D+}^\times=\{1\}$).

\noindent 
$\mathrm{(i)}$ Si  $F=\Q$, $\mathfrak{n}=p\Z$, avec $p$ un nombre premier,  
on a $p\!-\!1$ $(R,\mathfrak{n})$-pointes, au-dessus de la  $R$-pointe $\infty$
($\mathfrak{b}=\Z$), dont 

$\relbar$ $(p\!-\!1)/2$ non-ramifi{\'e}es, avec $\mathfrak{b'}=\Z$ et 
$\mathfrak{o}_{\overline{\mathcal{C}}}^\times=
\mathfrak{o}_{\mathcal{C}}^\times=\{1\}$. Chacune de ces pointes 
est seule dans sa $(R,\mathfrak{n})$-composante.

$\relbar$ $(p\!-\!1)/2$ ramifi{\'e}es, avec $\mathfrak{b'}=p^{-1}\Z$ et 
$\mathfrak{o}_{\overline{\mathcal{C}}}^\times=\{\pm 1\}\supset
\mathfrak{o}_{\mathcal{C}}^\times=\{1\}$, contenues dans une  seule
$(R,\mathfrak{n})$-composante.

\smallskip
\noindent
$\mathrm{(ii)}$ Si   $\mathfrak{n}=\mathfrak{p}^2$, avec $\mathfrak{p}$ un 
id{\'e}al   premier de $\mathfrak{o}$ de degr{\'e} r{\'e}siduel $1$ ($\N(\mathfrak{p})=p$, 
avec $p$ un nombre premier),  on a 3 types de    $(R,\mathfrak{n})$-pointes, 
au-dessus de la $R$-pointe $\infty$  ($\mathfrak{b}=\mathfrak{o}$) :

$\relbar$ si $\mathfrak{b'}=\mathfrak{o}$, on a $n=1$, 
$\mathfrak{o}_{\overline{\mathcal{C}}}^\times=
\mathfrak{o}_{\mathcal{C}}^\times=\mathfrak{o}^\times_{\mathfrak{p}^2}$, 
et donc  on a  $p(p\!-\!1)/[\mathfrak{o}^\times:
\mathfrak{o}^\times_{\mathfrak{p}^2}]$ pointes non-ramifi{\'e}es, 
chacune seule dans   sa $(R,\mathfrak{n})$-composante.

$\relbar$ si $\mathfrak{b'}=\mathfrak{p}^{-1}$, on a $n=p$, 
$\mathfrak{o}_{\overline{\mathcal{C}}}^\times=
\mathfrak{o}_{\mathcal{C}}^\times=\mathfrak{o}^\times_{\mathfrak{p}}$, et donc 
on a  $(p-\!1)^2/[\mathfrak{o}^\times:\mathfrak{o}^\times_{\mathfrak{p}}]$ 
pointes  peu ramifi{\'e}es, partag{\'e}es par groupes de $(p\!-\!1)$, en 
 $(p\!-\!1)/[\mathfrak{o}^\times: \mathfrak{o}^\times_{\mathfrak{p}}]$  
$(R,\mathfrak{n})$-composantes.

$\relbar$ si $\mathfrak{b'}=\mathfrak{p}^{-2}$, on a $n=p^2$, 
$\mathfrak{o}_{\overline{\mathcal{C}}}^\times=\mathfrak{o}^\times$, 
$\mathfrak{o}_{\mathcal{C}}^\times=\mathfrak{o}^\times_{\mathfrak{p}^2}$, 
et donc  on a  $p(p\!-\!1)/[\mathfrak{o}^\times:
\mathfrak{o}^\times_{\mathfrak{p}^2}]$ pointes tr{\`e}s
ramifi{\'e}es, contenus dans une  seule
$(R,\mathfrak{n})$-composante.

\smallskip
\noindent
$\mathrm{(iii)}$ Si   $\mathfrak{n}=\mathfrak{p}$, avec $\mathfrak{p}$ un 
id{\'e}al   premier de $\mathfrak{o}$ de degr{\'e} r{\'e}siduel $2$ ($\N(\mathfrak{p})=p^2$, 
avec $p$ un nombre premier), on a 2 types de 
  $(R,\mathfrak{n})$-pointes, au-dessus de la $R$-pointe $\infty$ 
($\mathfrak{b}=\mathfrak{o}$) :

$\relbar$ si $\mathfrak{b'}=\mathfrak{o}$, on a $n=1$, 
$\mathfrak{o}_{\overline{\mathcal{C}}}^\times=
\mathfrak{o}_{\mathcal{C}}^\times=\mathfrak{o}^\times_{\mathfrak{p}}$, et donc 
on a  $(p^2\!-\!1)/[\mathfrak{o}^\times:
\mathfrak{o}^\times_{\mathfrak{p}}]$ pointes non-ramifi{\'e}es, 
chacune seule dans   sa $(R,\mathfrak{n})$-composante.

$\relbar$ si $\mathfrak{b'}=\mathfrak{p}^{-1}$, on a $n=p$, 
$\mathfrak{o}_{\mathcal{C}}^\times=\mathfrak{o}^\times_{\mathfrak{p}}$
$\mathfrak{o}_{\overline{\mathcal{C}}}^\times=\{u\in \mathfrak{o}^\times
| u^p-u\in \mathfrak{p} \}$, et donc  
on a  $(p^2\!-\!1)/[\mathfrak{o}^\times:
\mathfrak{o}^\times_{\mathfrak{p}}]$ pointes 
peu ramifi{\'e}es, partag{\'e}es par groupes de 
$(p\!-\!1)/[\mathfrak{o}_{\overline{\mathcal{C}}}^\times:
\mathfrak{o}^\times_{\mathfrak{p}}]$, en  $(p+1)/[\mathfrak{o}^\times:
\mathfrak{o}_{\overline{\mathcal{C}}}^\times]$
$(R,\mathfrak{n})$-composantes.

\end{ex}

\smallskip

\section{Construction des cartes locales.} \label{carte-locale}

Soit $\mathfrak{c}$ un id{\'e}al de $F$, muni de sa
positivit{\'e} naturelle $\mathfrak{c}_+=\mathfrak{c}\cap(F\otimes\R)_+$. 
Posons  $\Delta=\Delta_F\N(\mathfrak{n})=\N(\mathfrak{dn})$.
Nous identifions $T_1\times D$ et $T$ par $(u,\epsilon)\mapsto 
\begin{pmatrix} u\epsilon & 0 \\ 0 & u^{-1} \end{pmatrix}$.

On on consid{\`e}re le  foncteur contravariant 
$\underline{\mathcal{M}}^1$  (resp. $\underline{\mathcal{M}}$)  de la
cat{\'e}gorie des   $\Z[{\frac{1}{\N(\mathfrak{n})}}]$-sch{\'e}mas  vers celle des ensembles,  
qui {\`a} un sch{\'e}ma $S$ associe 
l'ensemble des classes d'isomorphisme de quadruplets $(A,\iota,\lambda,\alpha)/S$ 
(resp. $(A,\iota,\overline{\lambda},\alpha)/S$), o{\`u} 
$(A,\iota)$ est une VAHB 
 (voir \cite{dimtildg} D{\'e}f.2.2), $\lambda$ est une $\mathfrak{c}$-polarisation 
 sur $A$ (resp. $\overline{\lambda}$ est une classe de
 $\mathfrak{c}$-polarisations; voir \cite{dimtildg} D{\'e}f.2.3), et 
 $\alpha: (\mathfrak{o}/\mathfrak{n})(1)\hookrightarrow A[\mathfrak{n}]$ 
est une  $\mu_{\mathfrak{n}}$-structure de niveau 
(voir \cite{dimtildg} D{\'e}f.2.5).

Le foncteur $\underline{\mathcal{M}}^1$  est 
repr{\'e}sentable par un $\Z[{\frac{1}{\N(\mathfrak{n})}}]$-sch{\'e}ma 
quasi-projectif, normal,   g{\'e}om{\'e}triquement connexe $M^1$ de 
dimension $d_F$, qui est  lisse 
au-dessus de   $\Z[{\frac{1}{\Delta}}]$ et 
muni d'un quadruplet universel 
$(\mathcal{A},\iota,\lambda,\alpha)$ (voir \cite{dimtildg} Thm.4.1).

Le foncteur $\underline{\mathcal{M}}$ admet un sch{\'e}ma de modules 
grossier $M$  sur $\Z[{\frac{1}{\N(\mathfrak{n})}}]$  
quasi-projectif,  normal, g{\'e}om{\'e}triquement connexe  
 et lisse au-dessus de   $\Z[{\frac{1}{\Delta}}]$ 
(voir \cite{dimtildg} Cor.4.2).

Le sch{\'e}ma $M$ est le quotient 
de $M^1$ par le groupe fini $\mathfrak{o}_{D+}^\times/
(\mathfrak{o}_{D+}^\times\cap\mathfrak{o}_{\mathfrak{n}}^{\times 2})$
qui agit proprement et librement par 
$[\epsilon]:(\mathcal{A},\iota,\lambda,\alpha)/S \mapsto 
(\mathcal{A},\iota,\epsilon\lambda,\alpha)/S$.
\vspace{5mm}

Le but de cette partie est de munir les  VAHB construites dans
la partie \ref{vahb-degenerantes} de diff{\'e}rentes 
$\mu_{\mathfrak{n}}$-structures de niveau, 
et ainsi fournir les cartes locales servant {\`a} compactifier la vari{\'e}t{\'e} 
modulaire de Hilbert $M$.

A chaque $(R,\mathfrak{n})$-composante  $\mathcal{C}$, on peut associer 
par la D{\'e}f.\ref{Rn-pointe} et la Prop.\ref{deuxideaux} des
id{\'e}aux $\mathfrak{b}$, $\mathfrak{b}'$ et  $X=
\mathfrak{cbb'}$, un entier $n$ {\'e}gal {\`a}  l'exposant
 du groupe $\mathfrak{b}'/\mathfrak{b}$,
des groupes d'unit{\'e}s $\mathfrak{o}_{\mathcal{C}}^\times$,
$\mathfrak{o}_{\overline{\mathcal{C}}}^\times$, 
$\mathfrak{o}_{\mathcal{C},1}^\times$,
$\mathfrak{o}_{\overline{\mathcal{C}},1}^\times$ et des sous-groupes 
$H_{\mathcal{C}}=\mathfrak{o}_{\overline{\mathcal{C}}}^\times/
\mathfrak{o}_{\mathcal{C}}^\times $,
$H_{\mathcal{C},1}=\mathfrak{o}_{\overline{\mathcal{C}},1}^\times/
\mathfrak{o}_{\mathcal{C},1}^\times $ du groupe   $(\Z/n\Z)^\times$
(ces objets sont {\it a priori} associ{\'e}s {\`a} une 
$(R,\mathfrak{n})$-pointe, mais sont constants au sein d'une 
$(R,\mathfrak{n})$-composante).

Soit une $(R,\mathfrak{n})$-composante  $\mathcal{C}$ et consid{\'e}rons le 
tore   $S=S_{\mathcal{C}}= \G_m\otimes X^*$.
Soit $\Sigma^{\mathcal{C}}$ un {\'e}ventail complet de $X_{+}^*$. 
Soit $\sigma\in\Sigma^{\mathcal{C}}$. 
La construction de la partie pr{\'e}c{\'e}dente, 
appliqu{\'e}e {\`a} $(X,\mathfrak{a},\mathfrak{b})$, nous
donne alors un sch{\'e}ma semi-ab{\'e}lien $G_{\sigma}/\overline{S}_{\sigma}$,
muni d'une action de  $\mathfrak{o}$ et  
dont la restriction {\`a} $G_{\sigma}^0/\overline{S}{}^{0}_{\sigma}$
est une VAHB $\mathfrak{c}$-polaris{\'e}e. 

En appliquant une deuxi{\`e}me fois la construction de la partie pr{\'e}c{\'e}dente,
cette fois {\`a}  $(X,\mathfrak{a},\mathfrak{b}')$, on obtient 
un sch{\'e}ma semi-ab{\'e}lien ${G_{\sigma}'}/\overline{S}_{\sigma}$,
muni d'une action de  $\mathfrak{o}$ et  
dont la restriction  $G'{}_{\!\!\sigma}^0/\overline{S}{}^{0}_{\sigma}$
est une VAHB $\mathfrak{c'}=\mathfrak{ab'}^{-1}$-polaris{\'e}e. Par fonctorialit{\'e} 
on a une fl{\`e}che $G_{\sigma}\rightarrow {G_{\sigma}'}$,
dont la restriction $G_{\sigma}^0\rightarrow G'{}_{\!\!\sigma}^0$
est une isog{\'e}nie.  On en  d{\'e}duit la suite exacte :
\begin{equation}\label{ntorsion2}
0 \rightarrow \mathfrak{b}'/\mathfrak{b} \overset{q}{\rightarrow} 
G_{\sigma}^0[\mathfrak{n}]\rightarrow 
G'{}_{\!\!\sigma}^0[\mathfrak{n}] \rightarrow 1.
\end{equation}

Consid{\'e}rons d'abord le cas o{\`u}  $\mathcal{C}$ est  non-ramifi{\'e}e. 
On a alors $\mathfrak{b}=\mathfrak{b}'$ et donc $X=\mathfrak{ab}$.
La vari{\'e}t{\'e} ab{\'e}lienne $G_{\sigma}^0$ associ{\'e}e {\`a} une 
$(R,\mathfrak{n})$-composante non-ramifi{\'e}e 
est naturellement munie d'une $\mu_{\mathfrak{n}}$-structure de niveau 
$(\mathfrak{o}/\mathfrak{n})(1) \cong
(\mathfrak{a}/\mathfrak{na})(1)\hookrightarrow G_{\sigma}^0[\mathfrak{n}]$, 
o{\`u} la premi{\`e}re fl{\`e}che vient de l'isomorphisme 
$\beta : \mathfrak{n}^{-1}\mathfrak{d}^{-1}/\mathfrak{d}^{-1}
\cong \mathfrak{n}^{-1}\mathfrak{a}^*/\mathfrak{a}^*$ et 
la deuxi{\`e}me du (\ref{ntorsion}).

Passons maintenant au cas o{\`u}  $\mathcal{C}$ est  ramifi{\'e}e. 
Afin de munir $G_{\sigma}^0$   d'une $\mu_{\mathfrak{n}}$-structure de niveau, on doit :

\smallskip
$\relbar$ choisir un rel{\`e}vement de $\mathfrak{b'}/\mathfrak{b} $
dans $\im(\beta)$ (appel{\'e} {\it uniformisation} de $\mathcal{C}$), 

\smallskip
$\relbar$ se placer dans ce cas au-dessus de 
$\Spec(\Z[\frac{1}{\N(\mathfrak{n})},\zeta_{\mathcal{C}}])$, 
o{\`u} $\zeta_{\mathcal{C}}$ d{\'e}signe  une racine de l'unit{\'e}  d'ordre {\'e}gal {\`a} 
l'exposant $n$ du groupe ab{\'e}lien  $\mathfrak{b'}/\mathfrak{b}$.
 
\medskip
Au-dessus de  $\Spec(\Z[\frac{1}{\N(\mathfrak{n})},\zeta_{\mathcal{C}}])$
on a un isomorphisme canonique $\mathfrak{b}^*/\mathfrak{b'}^*
\cong (\mathfrak{b'}/\mathfrak{b})(1)$, d'o{\`u} une 
$\mu_{\mathfrak{n}}$-structure de niveau sur $G_{\sigma}^0$ :

$ (\mathfrak{o}/\mathfrak{n})(1)   \hookrightarrow
(\mathfrak{a}/\mathfrak{na})(1)\times (\mathfrak{b}^*/\mathfrak{b'}^*)(1) 
\cong (\mathfrak{a}/\mathfrak{na})(1)\times \mathfrak{b'}/\mathfrak{b} 
\overset{(\ref{ntorsion})(\ref{ntorsion2})}
{\hookrightarrow}G_{\sigma}^0[\mathfrak{n}], $

\smallskip
\noindent  o{\`u} la premi{\`e}re inclusion vient de la fl{\`e}che 
$\beta:\mathfrak{n}^{-1}\mathfrak{d}^{-1}/\mathfrak{d}^{-1} \hookrightarrow
 \mathfrak{n}^{-1}\mathfrak{a}^*/\mathfrak{a}^*\times 
  \mathfrak{b'}/\mathfrak{b} $.

\vspace{-.4cm}
\begin{prop}\label{cartes-locales}

$\mathrm{(i)}$  Pour  toute $(R,\mathfrak{n})$-composante
 uniformis{\'e}e  $\mathcal{C}$ et pour tout
c{\^o}ne $\sigma\in\Sigma^\mathcal{C}$ la construction ci-dessus 
donne  un carr{\'e} cart{\'e}sien :

$$\begin{array}{cccc} G_{\sigma}^0\times 
\Spec(\Z[\frac{1}{\N(\mathfrak{n})},\zeta_{\mathcal{C}}])&\rightarrow & \mathcal{A} &\\
\downarrow&&\downarrow\\ \overline{S}{}^{0}_{\sigma}\times
\Spec(\Z[\frac{1}{\N(\mathfrak{n})},\zeta_{\mathcal{C}}]) 
&\rightarrow & M^1 &\rightarrow M\end{array}.$$

$\mathrm{(ii)}$ Changer l'uniformisation de la pointe $\mathcal{C}$
revient {\`a} se donner un {\'e}l{\'e}ment  
 $x\in (\mathfrak{ab})^*/(\mathfrak{ab'})^*=\Hom(\mathfrak{b'}/\mathfrak{b},
\mathfrak{n}^{-1}\mathfrak{a}^*/\mathfrak{a}^*)$ et correspond  donc {\`a} 
l'automorphisme de $\overline{S}{}^{0}_{\sigma}\times
\Spec(\Z[\frac{1}{\N(\mathfrak{n})},\zeta_{\mathcal{C}}])$ qui envoie
$q^\xi$ sur $\zeta_{\mathcal{C}}^{n\Tr_{F\!/\!\Q}(\xi x)} q^\xi$ 
($\xi\in \mathfrak{ab'}$).

\smallskip
$\mathrm{(iii)}$    Soient $\mathcal{C}_1$, $\mathcal{C}_2$ 
deux $(R,\mathfrak{n})$-composantes uniformis{\'e}es et 
soient deux c{\^o}nes $\sigma_i\subset X_{i,\R}^*$, $i=1,2$. 
Supposons qu'il  existe 

\smallskip
$\relbar$ un isomorphisme de 
$(R,\mathfrak{n})$-composantes $\mathcal{C}_1\cong \mathcal{C}_2$ 
(d'o{\`u}  $\xi\in F^\times$ tel que  $\mathfrak{a}_2^*=\xi \mathfrak{a}_1^*$, 
$\mathfrak{b}_2=\xi^{-1} \mathfrak{b}_1$ et  
$X_2^*=\xi^2 X_1^*$) induisant sur $\mathfrak{c}^*$ (via les polarisations de 
$L$ et $L'$), la multiplication par une  unit{\'e} 
$\epsilon\in \mathfrak{o}_{D+}^\times$,

\smallskip
$\relbar$   des {\'e}l{\'e}ments  $(u,\epsilon)\in\mathfrak{o}_{\mathcal{C}_1}^\times=
\mathfrak{o}_{\mathcal{C}_2}^\times$  et   $h\in H_{\mathcal{C}}$,
tels que  $\sigma_2=u^2 \epsilon \xi^2 \sigma_1$  et 
$\zeta_{\mathcal{C}_2}=\zeta_{\mathcal{C}_1}^h$.

\smallskip
Alors, on a un isomorphisme $\overline{S}{}^0_{\sigma_1}\times
\Spec(\Z[\frac{1}{\N(\mathfrak{n})},\zeta_{\mathcal{C}_1}])\cong 
\overline{S}{}^0_{\sigma_2}\times
\Spec(\Z[\frac{1}{\N(\mathfrak{n})},\zeta_{\mathcal{C}_2}])$ qui compl{\`e}te 
les  deux  fl{\`e}ches  $\overline{S}{}^0_{\sigma_i}\times
\Spec(\Z[\frac{1}{\N(\mathfrak{n})},\zeta_{\mathcal{C}_i}])
\rightarrow  M$  ($i=1,2$) du  $\mathrm{(i)}$  en 
un triangle commutatif.
\end{prop}

Le (i) et (ii) d{\'e}coulent de ce qui pr{\'e}c{\`e}de. Le (iii) utilise la fonctorialit{\'e} de 
la construction de $G_{\sigma}^0$ en $\sigma$ et  sa compatibilit{\'e} 
avec l'action de $\mathfrak{o}_{\mathcal{C}}^\times$ (voir
fin de la partie \ref{vahb-degenerantes} et la Prop.\ref{deuxideaux}(iv)).
\hfill$\square$

\medskip
Avant de d{\'e}crire la construction des compactifications toro{\"\i}dales 
arithm{\'e}tiques,  on doit la pr{\'e}parer. C'est l'objet des deux parties suivantes.

\section{Un th{\'e}or{\`e}me de descente formelle de Rapoport.}

La construction d'une compactification toro{\"\i}dale peut {\^e}tre vue comme l'ajout d'un bord
{\`a} $M$. On a un sch{\'e}ma formel de type fini candidat pour ce bord, 
{\`a} savoir l'analogue alg{\'e}brique de   :

$$\mathfrak{X}^{\an}=\coprod_{ \text{pointes } \mathcal{C} / \sim} 
\Big(\C^\times \otimes (\mathfrak{cbb}')^*
\Big)^\wedge_{\Sigma^\mathcal{C}}/\mathfrak{o}_{\mathcal{C}}^\times.$$
 
Le but de cette partie est de donner un crit{\`e}re abstrait, trouv{\'e}
par Rapoport \cite{rapoport},  pour r{\'e}soudre le probl{\`e}me
de ``{Descente Formelle}'', en l'occurrence, le probl{\`e}me
d'existence et unicit{\'e} du sch{\'e}ma recollement $Y$ 
d'un ouvert $Y^0$ et d'un sch{\'e}ma formel $\mathfrak{X}$ : 
$ Y^0\hookrightarrow Y\leftarrow \mathfrak{X}$.
Il repose en partie sur un crit{\`e}re d'immersion ouverte de Rapoport dont on rappellera l'{\'e}nonc{\'e}.

Le probl{\`e}me de Descente Formelle sera en fait d'abord pos{\'e} dans la cat{\'e}gorie des espaces alg{\'e}briques.
On verra dans la partie \ref{toroidale} 
 que les conditions d'application du crit{\`e}re sont satisfaites dans
notre cas.  

Dans cette partie $V$ d{\'e}signera un anneau de valuation discr{\`e}te complet, 
de corps des fractions $K$ et de corps r{\'e}siduel $k$. $S$ d{\'e}signe un $V$-sch{\'e}ma. 
\medskip

Soit $\Aff/S$ la cat{\'e}gorie des $S$-sch{\'e}mas affines, munie
de la topologie {\'e}tale. 
Un faisceau d'ensembles sur $\Aff/S$ s'appelle un $S$-espace.

\begin{defin} Une relation d'{\'e}quivalence {\'e}tale 
sur un $S$-sch{\'e}ma $U_1$ est donn{\'e}e par une  immersion ferm{\'e}e 
quasi-compacte $U_2\rightarrow U_1\times_S U_1$ de $S$-sch{\'e}mas
dont les deux projections sont {\'e}tales et qui d{\'e}finit une relation d'{\'e}quivalence :
pour tout $Y\in \Aff/S$, $U_2(Y)\rightarrow U_1(Y)\times_{S(Y)} 
U_1(Y)$ est une relation d'{\'e}quivalence. 

Un $S$-{\it espace alg{\'e}brique} est un $S$-espace qui est  quotient d'un sch{\'e}ma 
$U_1$, appel{\'e} un atlas {\'e}tale, par  une relation d'{\'e}quivalence {\'e}tale. 
\end{defin}

L'ensemble $\Alg/S$ des $S$-espaces alg{\'e}briques muni des 
fl{\`e}ches de $S$-espaces forme une cat{\'e}gorie.
On d{\'e}finit de m{\^e}me pour un sch{\'e}ma formel $S^\wedge$ 
la cat{\'e}gorie des $S^\wedge$-espaces  alg{\'e}briques formels, not{\'e}e $\Form/S^\wedge$.



\begin{defin} Soit $f:\mathfrak{X}'\rightarrow \mathfrak{X}$ un morphisme dans $\Form/S^\wedge$.
On dit que $f$ est un {\it {\'e}clatement admissible} de $\mathfrak{X}$ si $f$ 
est un {\'e}clatement $\mathfrak{X}'\rightarrow \mathfrak{X}$
 dans $\Form/S^\wedge$,  par rapport {\`a} un id{\'e}al 
qui contient une puissance de l'id{\'e}al de d{\'e}finition de $\mathfrak{X}$.

La cat{\'e}gorie des espaces rigides  $\Rig/S$ est la cat{\'e}gorie 
localis{\'e}e de $\Form/S$, par rapport aux {\'e}clatements admissibles.
\end{defin}



\begin{defin} Un {\it {\'e}paississement}
 de $(K,V)$ est un couple $(R,R^{(0)})$ tel que :

$\relbar$ $R$ est un anneau local artinien de corps r{\'e}siduel $K$. On note
$R_V$ l'image r{\'e}ciproque de $V$ dans $R$.

$\relbar$ $R^{(0)}\subset R_V \subset R$ est un sous-anneau noeth{\'e}rien tel que
le morphisme $R^{(0)}\rightarrow V$ soit surjectif et la localisation 
de $R^{(0)}$ au point g{\'e}n{\'e}rique de $V$ soit {\'e}gale {\`a} $R$ (c'est {\`a} dire
$R$ est le localis{\'e} de $R^{(0)}$ en $J=\ker(R^{(0)}\rightarrow V)$).
\end{defin}

Soit $\widetilde{\pi}$ un {\'e}l{\'e}ment de $R^{(0)}$ qui se projette sur
 une  uniformisante de $V$.
Pour tout $i\geq 1$ on pose 
$R^{(i)}=R^{(0)}\left[\frac{J}{\widetilde{\pi}^i}\right]$.
Alors $R^{(0)}\subset R^{(1)}\subset ...\subset R_V$ et $\cup_i R^{(i)}= R_V$.

On a $\Sprig(K)=\Spf(V)_{\rig}$ et $\Sprig(R):=\Spf(R^{(0)})_{\rig}$ 
($=\Spf(R^{(i)})_{\rig}$), car $\Spf(R^{(i)})$ est obtenu par
{\'e}clatement (admissible) de $\Spf(R^{(0)})$, par rapport {\`a} 
l'id{\'e}al $(\widetilde{\pi}^i)+J$ (car $J$ est nilpotent).

\begin{ex} 
 Soit l'anneau local artinien $R=K[t]/(t^2)$. 
Le sous-anneau  $R_V=V+K\cdot t$ n'est pas noeth{\'e}rien. Consid{\'e}rons
le sous-anneau  noeth{\'e}rien $R^{(0)}=V[t]/(t^2)$.
Alors $(R,R^{(0)})$ est un {\'e}paississement de $(K,V)$. 
On a $R^{(i)}=V+V\cdot \frac{t}{\widetilde{\pi}^i}$ et donc $\cup_i R^{(i)}= R_V$.
\end{ex}



A toute fl{\`e}che $f_{\rig}:\mathfrak{X}_{\rig}\rightarrow \mathfrak{Y}_{\rig}$
on peut associer un mod{\`e}le formel $f:\mathfrak{X}\rightarrow \mathfrak{Y}$,
d{\'e}fini {\`a} {\'e}clatement admissible pr{\`e}s. 

\begin{defin} 
$f_{\rig}$ est une {\it immersion ouverte}, s'il existe un mod{\`e}le formel 
$f$ qui est une immersion ouverte.
\end{defin} 

M. Rapoport a d{\'e}montr{\'e} le crit{\`e}re d'immersion ouverte suivant, qui est utilis{\'e} 
pour d{\'e}montrer le r{\'e}sultat de recollement abstrait que l'on a en vue.

\begin{theo}\label{immersion}{\rm(Th{\'e}or{\`e}me 3.15 de \cite{rapoport})}
$f_{\rig}$ est une immersion ouverte, si et seulement si, les 
deux conditions suivantes sont satisfaites :

\noindent $\mathrm{(i)_{rig}}$ 
Pour tout corps $K$, discr{\`e}tement valu{\'e}, l'application
$\Hom(\Sprig(K),\mathfrak{X}_{\rig})\overset{ f_{\rig*}}{\longrightarrow}
\Hom(\Sprig(K),\mathfrak{Y}_{\rig})$ est injective.

\noindent $\mathrm{(ii)_{rig}}$  
Pour tout {\'e}paississement $(R,R^{(0)})$ de $(K,V)$
on peut compl{\'e}ter de fa{\c c}on unique le diagramme commutatif suivant :
\xymatrix@R=15pt{ \Sprig(K)\ar[d]\ar[r] &\mathfrak{X}_{\rig}\ar[d] \\
\Sprig(R)\ar[r]\ar@{-->}[ru]  & \mathfrak{Y}_{\rig}}
\end{theo}

\begin{rque}
L'anneau $V$ {\'e}tant principal, il n'admet pas d'{\'e}clatements admissibles. 
La condition $\mathrm{(i)_{rig}}$ peut s'{\'e}crire donc 
$\Hom(\Spf(V),\mathfrak{X})\hookrightarrow
\Hom(\Spf(V),\mathfrak{Y})$, alors que le diagramme dans 
la condition $\mathrm{(ii)_{rig}}$ devient (pour $i$ assez grand) :

$$\xymatrix@R=15pt{ \Spf(V)\ar[d]\ar[rr] & &  \mathfrak{X}\ar[d] \\
\Spf(R^{(0)})  & \Spf(R^{(i)})\ar[l]\ar[r]\ar@{-->}[ru] & \mathfrak{Y}}$$
\end{rque}



Soit $S$ un sch{\'e}ma affine, de type fini sur le spectre d'un corps ou d'un anneau 
de Dedekind excellent (pour les applications aux compactifications 
toro{\"\i}dales, il suffit de prendre $S$ de type fini sur $\Z$).

Soit $A$ un anneau noeth{\'e}rien complet pour la topologie $I$-adique,
d{\'e}finie par un id{\'e}al $I\subset A$. Soit $\mathfrak{U}=\Spf(A)$
le sch{\'e}ma formel affine correspondant. Posons $\overline{U}=\Spec(A)$,
$\overline{U}_0=\Spec(A/I)$=l'{\it {\^a}me} de $\mathfrak{U}$ et 
$\overline{U}{}^0=\overline{U}\bs \overline{U}_0$.

\begin{lemme}\label{alg-form}{\rm(EGA III.5)}
Soit $Y$ un espace alg{\'e}brique de type fini 
sur $S$ et $Y_0\subset Y$ un sous-espace ferm{\'e}.
On suppose que $\overline{U}=\Spec(A)$ est un $S$-sch{\'e}ma et on se donne 
un  $S$-morphisme 
formel adique   $\mathfrak{f}:\mathfrak{U}\rightarrow Y|_{Y_0}$. 

Alors, il existe un unique morphisme $f:\overline{U}\rightarrow Y$
dont le compl{\'e}t{\'e} formel est  $\mathfrak{f}$.
\end{lemme}


\begin{defin} 
Un morphisme $g^0:\Spec(K)\rightarrow \overline{U}{}^0$ sera dit {\it permis},
s'il vient (via le lemme \ref{alg-form}) d'un morphisme formel de 
type fini $\mathfrak{g}:\Spf(V)\rightarrow \mathfrak{U}$.

\smallskip

Plus g{\'e}n{\'e}ralement (si $\overline{U}$ est un $S$-sch{\'e}ma), un morphisme 
$f^0:\overline{U}{}^0\rightarrow Y^0$ dans un espace alg{\'e}brique de type
fini sur $S$ sera dit {\it permis}, s'il existe une immersion ouverte
de $Y^0$ dans un $S$-espace alg{\'e}brique propre $Y$, telle que :
pour tout morphisme permis $\Spec(K)\rightarrow \overline{U}{}^0$,
l'unique extension {\`a} $\Spec(V)$ du morphisme compos{\'e} $\Spec(K)\rightarrow Y$,
envoie le point sp{\'e}cial dans $Y\bs Y^0$. 

\end{defin} 

Un  morphisme  $f^0$, provenant par restriction d'un morphisme 
$f:U\rightarrow Y$, est permis, s'il existe un morphisme formel
$\mathfrak{f}:\mathfrak{U}\rightarrow \mathfrak{Y}=Y|_{Y_0}$ qui fait
commuter le diagramme suivant :
\xymatrix@R=10pt{ \mathfrak{U}\ar[d]\ar^{\mathfrak{f}}[r]& \mathfrak{Y}\ar[d]\\ 
\overline{U}\ar[d]\ar^{f}[r] &Y\ar[d] \\
\overline{U}{}^0 \ar^{f^0}[r] & Y^0}

En d'autres termes, un morphisme est permis 
s'il ``envoie le bord sur le bord''.


\begin{defin} \label{decoupage}
Soit  $\mathfrak{X}$ un $S$-espace alg{\'e}brique formel, s{\'e}par{\'e}
et de type fini.
 Un {\it d{\'e}coupage} de $\mathfrak{X}$ est la donn{\'e}e :

$\relbar$d'un atlas affine  $\mathfrak{U}_2=\Spf(A_2)
\rightrightarrows \mathfrak{U}_1=\Spf(A_1)  \rightarrow \mathfrak{X}$, et 

$\relbar$d'un espace alg{\'e}brique $Y^0$ de type fini sur $S$, tel que 
les deux compos{\'e}s suivants soient {\'e}gaux : $\overline{U}{}_2^0\rightrightarrows
 \overline{U}{}_1^0 \overset{f^0}{\rightarrow}  Y^0$, o{\`u}  
$\overline{U}_1=\Spec(A_1)$ et $\overline{U}_2=\Spec(A_2)$ et  les fl{\`e}ches   
$\overline{U}_2 \rightrightarrows\overline{U}_1$ viennent, via  le lemme 
\ref{alg-form}, des fl{\`e}ches $\mathfrak{U}_2\rightrightarrows \mathfrak{U}_1$.

Le d{\'e}coupage est dit {\it effectif}, s'il existe un $S$-espace alg{\'e}brique 
de type fini $Y$, une immersion ouverte 
$j:Y^0\hookrightarrow Y$ et un isomorphisme 
$\varphi:\mathfrak{X} \overset{\sim}{\rightarrow} \mathfrak{Y}$,
o{\`u}  $\mathfrak{Y}$ est 
le compl{\'e}t{\'e} formel de $Y$ le long de $Y\bs Y^0$, tels que le 
morphisme $f:\overline{U}_1\rightarrow Y$, venant (via le lemme \ref{alg-form})
du morphisme $\mathfrak{f}:\mathfrak{U}_1\rightarrow \mathfrak{X} 
\overset{\sim}{\rightarrow} \mathfrak{Y}$, induise 
$f^0:\overline{U}{}_1^0\rightarrow Y^0$ sur $\overline{U}{}_1^0\subset \overline{U}_1$.
$$\xymatrix@R=15pt{ \mathfrak{U}_2\ar@<0.5ex>[r] \ar@<-0.5ex>[r]\ar[d] & \mathfrak{U}_1
\ar^{\mathfrak{f}}[r]\ar[d] &  \mathfrak{X}\cong  \mathfrak{Y} \ar[d]\\
\overline{U}_2\ar@<0.5ex>[r] \ar@<-0.5ex>[r] &\overline{U}_1 \ar^{f}[r] & Y \\
\overline{U}{}_2^0\ar@<0.5ex>[r]\ar@{^{(}->}[u] \ar@<-0.5ex>[r] &
\overline{U}{}_1^0 \ar^{f^0}[r]\ar@{^{(}->}[u]  &   Y^0 \ar@{^{(}->}[u]}$$
\end{defin}


\begin{theo}{\rm(Th{\'e}or{\`e}me 3.5 de \cite{rapoport})}\label{algebrisation}
Soit un d{\'e}coupage. On suppose :

$\relbar$  $\overline{U}{}_1^0$ est sch{\'e}matiquement dense
 dans $\overline{U}_1$ (i.e. $\mathcal{O}_{\overline{U}_1}\hookrightarrow 
\mathcal{O}_{\overline{U}{}_1^0}$).

$\relbar$ $Y^0$ est compactifiable (i.e. il existe une 
$S$-immersion ouverte $Y^0\hookrightarrow Y^*$ avec  $Y^*$ propre sur $S$).

$\relbar$ Le morphisme $f^0 : \overline{U}{}_1^0  \rightarrow Y^0$ 
est permis.

$\relbar$ Pour tout anneau de valuation discr{\`e}te complet $V$, de corps 
des fractions $K$:

$\mathrm{(i')_{rig}}$ la suite 
$\overline{U}{}_2^0(K)_{\mathrm{permis}} \rightrightarrows
\overline{U}{}_1^0(K)_{\mathrm{permis}} 
\rightarrow Y^0(K)_{\mathrm{permis}}$ est exacte, et 

$\mathrm{(ii')_{rig}}$  pour tout {\'e}paississement $(R,R^{(0)})$   de $(K,V)$
on peut compl{\'e}ter de fa{\c c}on unique le diagramme commutatif suivant :
\xymatrix@R=15pt{ \Spec(K)\ar[d]\ar^{\hspace{3mm}\mathrm{permis}}[r] 
& \overline{U}{}_1^0 \ar[d] \\
\Spec(R)\ar^{\hspace{3mm}\mathrm{permis}}[r]\ar@{-->}[ru]  & Y^0}.

Alors le d{\'e}coupage est effectif.
\end{theo}

\vspace{-.4cm}
\section{La construction de Raynaud.}\label{raynaud2}

Pour pouvoir v{\'e}rifier les conditions 
$\mathrm{(i')_{rig}}$ et $\mathrm{(ii')_{rig}}$ ci-dessus
dans la situation o{\`u} l'ouvert $Y^0$ est l'espace de modules $M^1$ et le sch{\'e}ma  formel
$\mathfrak{X}$ est celui donn{\'e} par les cartes locales de la proposition 
\ref{cartes-locales},
on a besoin de la construction suivante (donn{\'e}e par Raynaud dans 
\cite{raynaud} et reprise par Rapoport dans le cas d'une VAHB \cite{rapoport}). 
Il est {\`a} noter que l'on a  besoin de cette construction non 
seulement sur un corps mais aussi  sur un {\'e}paississement artinien, auquel cas
l'argument donn{\'e} par  Raynaud reste valable.

Soit $V$ un anneau de valuation discr{\`e}te complet de corps
des fractions $K$, et soit $(R,R^{(0)})$ un {\'e}paississement de $(K,V)$.

\begin{defin}
Une vari{\'e}t{\'e} ab{\'e}lienne  
$A$ sur $R$ (resp. sur $K$) est dite {\`a} r{\'e}duction semi-stable
(d{\'e}ploy{\'e}e) s'il existe un sch{\'e}ma en groupes lisse sur 
 $R^{(i)}$, pour un certain $i\geq 0$, (resp. sur $V$), prolongeant
$A$ et dont la fibre sp{\'e}ciale est une extension d'une vari{\'e}t{\'e} 
ab{\'e}lienne par un tore (d{\'e}ploy{\'e}). 
\end{defin}

Pour des raisons de dimension, si une VAHB 
sur $R$ (ou sur $K$) est {\`a}  mauvaise r{\'e}duction semi-stable d{\'e}ploy{\'e}e,
alors la fibre sp{\'e}ciale est un tore d{\'e}ploy{\'e}. Dans ce cas 
la description rigide-analytique de Raynaud devient :

\begin{theo}$\mathrm{(Raynaud)}$
Soit A une VAHB sur $R$ (ou sur $K$) {\`a} mauvaise r{\'e}duction semi-stable
d{\'e}ploy{\'e}e. Alors :
$$A_{\rig}=(\Gm\otimes \mathfrak{a}^*)_{\rig}/\mathfrak{b}_{\rig},$$
o{\`u} $\mathfrak{a}$ et $\mathfrak{b}$  sont des id{\'e}aux de $F$. De plus :

$\relbar$ on a une suite exacte 
$0 \rightarrow (\mathfrak{a}/\mathfrak{na})(1) \rightarrow A[\mathfrak{n}]
\rightarrow  \mathfrak{n}^{-1}\mathfrak{b}/\mathfrak{b}\rightarrow 0$,

$\relbar$ la forme bilin{\'e}aire $\langle\enspace,\enspace\rangle:
\mathfrak{a}\times \mathfrak{b}\rightarrow 
\mathrm{Val(K)}\cong\Z$ \enspace $(\alpha,\beta)\mapsto 
\mathrm{val}(\mathfrak{X}^\alpha(\beta))$ v{\'e}rifie
$\langle a\alpha,\beta\rangle=\langle \alpha,a\beta\rangle$
 pour tout $a\in\mathfrak{o}$,
et donc d{\'e}finit un unique {\'e}l{\'e}ment $\xi^*\in (\mathfrak{ab})^*$,
{\`a} $\Q_+^\times$ pr{\`e}s et {\`a} l'action de $\mathfrak{o}^{\times}$ pr{\`e}s,

$\relbar$ le c{\^o}ne positif des polarisations sur $A$,
$\mathcal{P}(A)\subset \Sym_{\mathfrak{o}}(A,A^t)=\mathfrak{ab}^{-1}$
est obtenu comme produit de l'unique positivit{\'e} sur
$\mathfrak{ab}$ pour laquelle $\xi^*>0$ et de la positivit{\'e} naturelle
sur $\mathfrak{b}^{-2}$.
\end{theo}

\vspace{-.6cm}
\section{Compactifications toro{\"\i}dales arithm{\'e}tiques.}\label{toroidale}

\vspace{4mm}

\subsubsection*{Construction  des compactifications toro{\"\i}dales.}

\vspace{-6mm}
\begin{defin} Un {\it {\'e}ventail $\Gamma$-admissible} 
$\Sigma=(\Sigma^{\mathcal{C}})_{\mathcal{C}}$ 
est la donn{\'e}e pour chaque $(R,\mathfrak{n})$-composante  
$\mathcal{C}$  d'un {\'e}ventail complet
$\Sigma^{\mathcal{C}}$ de $X_{+}^*$, 
stable par $\mathfrak{o}_{\mathcal{C}}^\times$ et 
contenant un nombre fini d'{\'e}l{\'e}ments modulo
cette action, de sorte que les donn{\'e}es soient compatibles aux isomorphismes 
de $(R,\mathfrak{n})$-composantes $\mathcal{C}\cong \mathcal{C}'$.
\end{defin}

Voici l'analogue du  r{\'e}sultat principal de l'article \cite{rapoport} 
dans le cas de groupe de niveau $\Gamma$ (on rappelle que
 $\Gamma$ est sans torsion).

\begin{theo} 
Soit $\Sigma=\{\Sigma^{\mathcal{C}}\}_{\mathcal{C}}$ 
un {\'e}ventail  $\Gamma$-admissible.

$\mathrm{(i)}$ Il existe une immersion ouverte $j:M^1\hookrightarrow 
M^1_{\Sigma}$ de $\Spec(\Z[\frac{1}{\N(\mathfrak{n})}])$-sch{\'e}mas
et un isomorphisme de sch{\'e}mas formels 
$$\varphi: \coprod_{(R,\mathfrak{n})\mathrm{-composantes}/\sim} 
\left(S_{\Sigma^{\mathcal{C}}}^{\wedge}/\mathfrak{o}_{\mathcal{C},1}^{\times} \right)\times \textstyle
\Spec(\Z[\frac{1}{\N(\mathfrak{n})},\zeta_{\mathcal{C}}]^{H_{\mathcal{C},1}})
\overset{\sim}{\longrightarrow}  M_{\Sigma}^{1\wedge},$$
(o{\`u} $M_{\Sigma}^{1\wedge}$ est le compl{\'e}t{\'e} formel
de $M^1_{\Sigma}$ le long de sa partie {\`a} l'infini), tels que 
pour toute $(R,\mathfrak{n})$-composante $\mathcal{C}$ et 
pour tout  $\sigma\in\Sigma^{\mathcal{C}} $ on a la propri{\'e}t{\'e} suivante:
l'image r{\'e}ciproque  de la VAHB universelle sur $M^1$ par le morphisme 
$\overline{S}_{\sigma}\times
 \Spec(\Z[\frac{1}{\N(\mathfrak{n})},\zeta_{\mathcal{C}}])
\rightarrow M^1_{\Sigma}$
(d{\'e}duit par le lemme \ref{alg-form} du morphisme formel 
$S_{\sigma}^{\wedge}\times
 \Spec(\Z[\frac{1}{\N(\mathfrak{n})},\zeta_{\mathcal{C}}])
\rightarrow  M_{\Sigma}^{1\wedge}$ construit 
{\`a} l'aide de $\varphi$), soit la VAHB $\mathfrak{c}$-polaris{\'e}e avec 
$\mu_{\mathfrak{n}}$-structure de niveau $G_{\sigma}^0\times
 \Spec(\Z[\frac{1}{\N(\mathfrak{n})},\zeta_{\mathcal{C}}])$   sur
$\overline{S}{}^{0}_{\sigma}\times
 \Spec(\Z[\frac{1}{\N(\mathfrak{n})},\zeta_{\mathcal{C}}])$ 
construite  dans la proposition \ref{cartes-locales}(i).
Le couple $(j,\varphi)$ est unique, {\`a} unique isomorphisme pr{\`e}s.

(ii) Il existe  une immersion ouverte $j:M\hookrightarrow 
M_{\Sigma}$ de $\Spec(\Z[\frac{1}{\N(\mathfrak{n})}])$-sch{\'e}mas
et un isomorphisme de sch{\'e}mas formels 
$$\varphi: \coprod_{(R,\mathfrak{n})\mathrm{-composantes}/\sim} 
\left(S_{\Sigma^{\mathcal{C}}}^{\wedge}/\mathfrak{o}_{\mathcal{C}}^{\times} \right)\times \textstyle
\Spec(\Z[\frac{1}{\N(\mathfrak{n})},\zeta_{\mathcal{C}}]^{H_{\mathcal{C}}})
\overset{\sim}{\longrightarrow}  M_{\Sigma}^{\wedge}.$$

\end{theo}

\noindent {\bf D{\'e}monstration : } (i)
Il y a un nombre fini de $(R,\mathfrak{n})$-composantes $\mathcal{C}$ 
modulo isomorphisme.
Soit $\{\sigma_i^\mathcal{C}\}$ un ensemble fini de repr{\'e}sentants des c{\^o}nes
de l'{\'e}ventail $\Sigma^{\mathcal{C}}$, modulo  
l'action de $\mathfrak{o}_{\mathcal{C},1}^{\times}$.

Consid{\'e}rons le sch{\'e}ma  formel affine 
 $\mathfrak{U}_1:=\underset{\mathcal{C}/\sim}{\coprod}
\underset{i}{\coprod} S_{\sigma_i^\mathcal{C}}^{\wedge}
\times \Spec(\Z[\frac{1}{\N(\mathfrak{n})},\zeta_\mathcal{C}])$.
Il est de type fini sur $\Z$ et muni d'un 
morphisme {\'e}tale (``immersions toro{\"\i}dales'' et quotient 
{\'e}tale par le groupe $H_{\mathcal{C},1}$) dans 
$\mathfrak{X}:=\underset{\mathcal{C}/\sim}{\coprod}
\left(S_{\Sigma^{\mathcal{C}}}^{\wedge}/
\mathfrak{o}_{\mathcal{C},1}^{\times} \right)
\times \Spec(\Z[\frac{1}{\N(\mathfrak{n})},\zeta_\mathcal{C}]^{H_{\mathcal{C},1}})$.

Posons $\overline{U}_1=\underset{\mathcal{C}/\sim}{\coprod}
\underset{i}{\coprod}\overline{S}_{\sigma_i^\mathcal{C}}
\times \Spec(\Z[\frac{1}{\N(\mathfrak{n})},\zeta_\mathcal{C}])$.

D'apr{\`e}s la proposition \ref{cartes-locales}(i) on a un  morphisme  
$f^0:\overline{U}{}_1^0 \rightarrow M^1$, qui  est permis, 
 car toute vari{\'e}t{\'e} ab{\'e}lienne  obtenue comme image 
r{\'e}ciproque, par un morphisme permis $\Spec(K)\rightarrow 
\overline{S}{}_{\sigma_i^\mathcal{C}}^0\times 
\Spec(\Z[\frac{1}{\N(\mathfrak{n})},\zeta_\mathcal{C}])$,
de la vari{\'e}t{\'e} ab{\'e}lienne  $G_{\sigma}^0\times 
\Spec(\Z[\frac{1}{\N(\mathfrak{n})},\zeta_{\mathcal{C}}])$
est {\`a} mauvaise r{\'e}duction d'apr{\`e}s la partie  \ref{vahb-degenerantes}.

Posons  $\mathfrak{U}_2:=\mathfrak{U}_1\times_\mathfrak{X}\mathfrak{U}_1=
\Spf(A_2)$ et  $\overline{U}_2=\Spec(A_2)$. 
Les deux fl{\`e}ches compos{\'e}es $\overline{U}{}_2^0\rightrightarrows 
\overline{U}{}_1^0 \rightarrow M^1$ sont {\'e}gales  par
compatibilit{\'e} de  la construction de Mumford avec les inclusions 
$\sigma'\subset\sigma$, avec l'action de $\mathfrak{o}_{\mathcal{C},1}^\times $
 et avec l'action de $H_{\mathcal{C},1}$
(appliquer la proposition \ref{cartes-locales}(iii) dans le cas $D=\Gm$).

\medskip

V{\'e}rifions  la condition  $\mathrm{(i')_{rig}}$ du th{\'e}or{\`e}me  
\ref{algebrisation} :

Soient $g_1^0,g_2^0:\Spec(K)\rightarrow \overline{U}{}_1^0$ deux morphismes
permis avec $f^0\circ g_1^0=f^0\circ g_2^0$.

Chaque morphisme $g_j^0$ se factorise par un certain 
$\overline{S}{}^0_{\sigma_j}\times 
\Spec(\Z[\frac{1}{\N(\mathfrak{n})},\zeta_{\mathcal{C}_j}])$,
o{\`u} $\sigma_j$ d{\'e}signe un des $\sigma_i^{\mathcal{C}_j}$ et  d{\'e}termine ainsi :

\noindent$\relbar$une $(R,\mathfrak{n})$-composante  $\mathcal{C}_j$ 
({\`a} laquelle sont attach{\'e}s des objets  
$\mathfrak{a}_j$, $\mathfrak{b}_j$, $\mathfrak{b}'_j$, 
$X_j$, $\beta_j$), 

\noindent$\relbar$une racine de l'unit{\'e} $\zeta_\mathcal{C}^{(j)}\in K$,
d'ordre l'exposant $n_j$ du groupe  $\mathfrak{b}'_j/\mathfrak{b}_j$,

\noindent$\relbar$un c{\^o}ne $\sigma_j$ de $\Sigma^{\mathcal{C}_j}$ et 
un morphisme $\psi_j:R_{\sigma_j}^{\wedge}\rightarrow V$,
d'o{\`u}  un {\'e}l{\'e}ment $\xi_j^*\in\sigma_j\cap X_j^*$, d{\'e}termin{\'e} par
la propri{\'e}t{\'e} suivante :
pour tout $\xi\in  \check{\sigma}_j\cap X_j$ 
on a $\mathrm{val}(\psi_j(q^{\xi}))=\Tr_{F\!/\!\Q}(\xi \xi_j^*)$.

\smallskip
Le morphisme permis  $f^0\circ g_{j}^0$ fournit 
une  VAHB $A$ sur $K$ munie d'une  $\mathfrak{c}$-polarisation et 
 $\mu_{\mathfrak{n}}$-structure de niveau, {\`a} mauvaise r{\'e}duction semi-stable
d{\'e}ploy{\'e}e.

L'uniformisation de Raynaud-Tate de la VAHB $A$, d{\'e}crite dans 
la partie \ref{raynaud2},    donne alors :

$\relbar$ deux  id{\'e}aux $\mathfrak{a}$ et $\mathfrak{b}$, tels que 
$A_{\rig}=(\Gm\otimes \mathfrak{a}^*)_{\rig}/\mathfrak{b}_{\rig}$ et 
 $\mathfrak{c}=\Sym_{\mathfrak{o}}(A,A^t)=\mathfrak{ab}^{-1}$
(ceci nous donne  une $R$-pointe  $\mathcal{C}$, bien d{\'e}finie modulo
isomorphisme). Comme la construction de Mumford  et celle de Raynaud 
sont inverses l'une {\`a} l'autre (i.e. le $1$-motif associ{\'e} par
 Raynaud {\`a} la  VAHB sur $K$ construite par Mumford est 
le $1$-motif du d{\'e}part), les $R$-pointes sous-jacentes 
de  $\mathcal{C}_1$ et $\mathcal{C}_2$ sont  isomorphes {\`a} $\mathcal{C}$.

$\relbar$ une suite exacte : 
$0 \rightarrow (\mathfrak{a}/\mathfrak{na})(1) 
\rightarrow A[\mathfrak{n}]\rightarrow 
\mathfrak{n}^{-1}\mathfrak{b}/\mathfrak{b}\rightarrow 0$. 
Ainsi, la  $\mu_{\mathfrak{n}}$-structure 
de niveau sur $A$ d{\'e}termine-t-elle bien une 
une $(R,\mathfrak{n})$-composante  $\mathcal{C}$ au-dessus de la 
$R$-pointe $\mathcal{C}$ et une racine de l'unit{\'e} $\zeta_\mathcal{C}$. 
De nouveau par compatibilit{\'e} de
la construction de Mumford  et celle de Raynaud on d{\'e}duit que 
$\zeta_{\mathcal{C}_1}$ et $\zeta_{\mathcal{C}_2}$ sont
conjugu{\'e}es sous $H_{\mathcal{C},1}$ et que 
les $(R,\mathfrak{n})$-composantes $\mathcal{C}_1$ et 
$\mathcal{C}_2$ sont isomorphes, et donc {\'e}gales,
car elles vivent dans un ensemble de repr{\'e}sentants modulo isomorphisme.

$\relbar$ un   {\'e}l{\'e}ment $\xi^*\in (\mathfrak{ab})^*_+ $
bien d{\'e}fini modulo   $\mathfrak{o}_{\mathcal{C},1}^\times$.
Un derni{\`e}re fois par 
 compatibilit{\'e} des constructions de Mumford  et  de Raynaud,
on trouve  que  $\xi_1^*\in\sigma_1$ et $\xi_2^*\in\sigma_2$ sont dans la 
m{\^e}me $\mathfrak{o}_{\mathcal{C},1}^\times$-orbite. Par cons{\'e}quent 
$\xi_1^*=\xi_2^*$ et, par exemple $\sigma_1\subset \sigma_2$.

On en d{\'e}duit qu'il existe 
un morphisme permis $h^0: \Spec(K)\rightarrow \overline{U}{}_2^0$
tel que $g_1^0=p_1^0\circ h^0$ et  $g_2^0=p_2^0\circ h^0$, ce
qui termine la v{\'e}rification  du  $\mathrm{(i')_{rig}}$.

\medskip

V{\'e}rifions  la condition  $\mathrm{(ii')_{rig}}$
du th{\'e}or{\`e}me  \ref{algebrisation} :

Les morphismes permis $\Spec(K)\rightarrow \overline{U}{}_1^0$
et $\Spec(R)\rightarrow M^1$ nous donnent deux VAHB  
$A/K$ et $A'/R$ {\`a} mauvaise r{\'e}duction, avec $A\cong A'\times_R K$. 
Comme dans la v{\'e}rification de $\mathrm{(i')_{rig}}$,
la fl{\`e}che $\Spec(K)\rightarrow \overline{U}{}_1^0$
d{\'e}termine les  des donn{\'e}s combinatoires  $\mathcal{C}=
(\mathfrak{a},\mathfrak{b},X,\beta), \zeta_{\mathcal{C}},\xi^*\in X^*$.
Par la th{\'e}orie de Raynaud-Tate $A$ et $A'$  admettent des uniformisation 
rigides analytiques 
$A_{\rig}=(\Gm\otimes \mathfrak{a}^*)_{\rig}/\mathfrak{b}_{\rig}$
(compatibilit{\'e} entre la construction de  Mumford et celle de Raynaud)
et  $A'_{\rig}=(\Gm\otimes {\mathfrak{a}'}^*)_{\rig}/\mathfrak{b}'_{\rig}$.
Comme $A_{\rig}=A'_{\rig}\times_R K$, on en d{\'e}duit que l'on peut
prendre $\mathfrak{a}=\mathfrak{a}'$, $\mathfrak{b}=\mathfrak{b}'$,
$\zeta_{\mathcal{C}}=\zeta_{\mathcal{C}'}$ et $\xi^*={\xi'}^*$,
d'o{\`u} le $\mathrm{(ii')_{rig}}$.

Nous sommes maintenant en mesure d'appliquer le th{\'e}or{\`e}me  \ref{algebrisation}
qui nous donne le couple cherch{\'e} $(j,\varphi)$, dont on admet l'unicit{\'e}.

(ii) Comme  $\Sigma^{\mathcal{C}}$ est 
stable par $\mathfrak{o}_{\mathcal{C}}^\times$ (et non-seulement par 
$\mathfrak{o}_{\mathcal{C},1}^\times$), le  groupe fini 
 $\mathfrak{o}_{D+}^\times/
(\mathfrak{o}_{D+}^\times\cap\mathfrak{o}_{\mathfrak{n}}^{\times 2})$
du rev{\^e}tement galoisien {\'e}tale
$M^1\rightarrow M$ agit proprement et librement  sur $M^1_{\Sigma}$ et 
la construction du (i) passe au quotient. La fl{\`e}che 
$M^1_{\Sigma}\rightarrow M_{\Sigma}$ est encore {\'e}tale.
\hfill $\square$

\begin{rque}
Soit $\Sigma=(\Sigma^{\mathfrak{b}})_{\mathfrak{b}\in \mathcal{I}_F}$, 
o{\`u} pour tout id{\'e}al $\mathfrak{b}$, $\Sigma^{\mathfrak{b}}$
est un {\'e}ventail $\mathfrak{o}^\times$-admissible de 
$(\mathfrak{cb}^2)_+^*$. On aurait pu tenter de d{\'e}finir 
$M_{\Sigma}$ comme la normalisation dans $M$ de la compactification de 
Rapoport $M(\mathfrak{c})_{\Sigma}$ de l'espace de modules $M(\mathfrak{c})$.
Le probl{\`e}me est que le sch{\'e}ma $M_{\Sigma}$ ainsi d{\'e}fini n'est jamais lisse.
En effet, pour compactifier chaque  $(R,\mathfrak{n})$-pointe $\mathcal{C}$ 
qui est au-dessus de la $R$-pointe correspondant {\`a}  $\mathfrak{b}$ 
on 
utilise le m{\^e}me {\'e}ventail $\Sigma^{\mathfrak{b}}$. 
Or, si $\mathfrak{bb'}^{-1}\neq n\mathfrak{o}$ ($n\in \Z$), 
$\Sigma^{\mathfrak{b}}$ ne peut pas {\^e}tre un {\'e}ventail lisse 
pour $(\mathfrak{cb}^2)_+^*$ et $(\mathfrak{cbb'})_+^*$ simultan{\'e}ment.
\end{rque}

\subsubsection*{Propri{\'e}t{\'e}s des compactifications toro{\"\i}dales.}

Dans la suite, pour all{\'e}ger les notations, nous {\'e}crirons 
$\overline{M}$ {\`a} la place de $M_\Sigma$, en gardant en t{\^e}te la
d{\'e}pendance du syst{\`e}me d'{\'e}ventails $\Sigma$.

\begin{cor}  Localement pour la topologie {\'e}tale sur 
$\Spec(\Z[\frac{1}{\N(\mathfrak{n})}])$,
 $j : M\hookrightarrow  \overline{M}$ est isomorphe {\`a} 
$S_{\mathcal{C}}\hookrightarrow  S_{\sigma}$ pour 
un certain couple $\mathcal{C}$, $\sigma\in \Sigma^{\mathcal{C}}$. 

\noindent En particulier, pour tout c{\^o}ne 
$\sigma\in \Sigma^{\mathcal{C}}\bs\{0\}$,  et tout corps alg{\'e}briquement clos $k$ 
de caract{\'e}ristique  $p$ ne divisant pas $\N(\mathfrak{n})$, 
l'ensemble des $k$-points de la 
strate $ M(\sigma)$ de $\overline{M}$ s'identifie {\`a} celui des $k$-points 
de la strate ferm{\'e}e $S(\sigma)$ de l'immersion torique affine 
$S\hookrightarrow S_\sigma$.
\end{cor}

Ceci r{\'e}sulte du fait 
que $\mathfrak{o}_{\mathcal{C}}^\times$ op{\`e}re librement sur l'ensemble des
strates non-ouvertes de  $S_{\mathcal{C}}\hookrightarrow 
S_{\Sigma^{\mathcal{C}}}$, et donc localement pour 
la topologie {\'e}tale 
$S_{\Sigma^{\mathcal{C}}}^\wedge/\mathfrak{o}_{\mathcal{C}}^\times $
(et donc $\overline{M}^\wedge$) est isomorphe {\`a} $S_\sigma^\wedge$, pour 
un certain $\sigma\in \Sigma^{\mathcal{C}}$.

\begin{cor}
 Quitte {\`a} raffiner l'{\'e}ventail $\Sigma$,  on  obtient un sch{\'e}ma 
 $\overline{M}$  qui est lisse au-dessus de   $\Spec(\Z[\frac{1}{\Delta}])$.
\end{cor}

\begin{prop}
 Il existe un unique sch{\'e}ma en groupes semi-ab{\'e}lien 
$\overline{f}:\mathfrak{G}\rightarrow \overline{M^1}$ qui prolonge 
la VAHB  universelle  $f:\mathcal{A}\rightarrow M^1$. 
Ce sch{\'e}ma en groupes est muni d'une action de 
$\mathfrak{o}$ et c'est un tore au-dessus de $\overline{M^1} \bs M^1$.
\end{prop}

\noindent {\bf D{\'e}monstration : } L'unicit{\'e} est montr{\'e}e dans un 
cadre beaucoup plus g{\'e}n{\'e}ral dans le chapitre I du livre de
Chai et Faltings \cite{FaCh}. 
Pour l'existence on consid{\`e}re le diagramme suivant :

$$\xymatrix@R=10pt@C=10pt{ 
& \mathcal{A}\ar@{-->}[rr]\ar[dd]& &\mathfrak{G}\ar@{-->}[dd] & & \\
G_{\sigma}^0[\frac{1}{\N(\mathfrak{n})},\zeta_{\mathcal{C}}] \ar[ru]\ar[rr]\ar[dd]& 
& G_{\sigma}[\frac{1}{\N(\mathfrak{n})},\zeta_{\mathcal{C}}]\ar[dd]\ar@{-->}[ru] & & & \\
 & M^1 \ar[rr] & & \overline{M^1} & &   \overline{M^1}^{\wedge}\ar[ll] \\
 \overline{S}{}^{0}_{\sigma}[\frac{1}{\N(\mathfrak{n})},\zeta_{\mathcal{C}}]\ar[rr]\ar[ru] & & 
\overline{S}_{\sigma}[\frac{1}{\N(\mathfrak{n})},\zeta_{\mathcal{C}}]\ar[ru] & & 
S_{\sigma}^{\wedge}[\frac{1}{\N(\mathfrak{n})},\zeta_{\mathcal{C}}]\ar[ru]\ar[ll]  & }$$

\begin{theo}
Le $\Spec(\Z[\frac{1}{\N(\mathfrak{n})}])$-sch{\'e}ma $\overline{M}$ est  propre.
\end{theo}

\noindent {\bf D{\'e}monstration : } L'id{\'e}e, comme dans \cite{rapoport},
 est d'appliquer  le crit{\`e}re valuatif de propret{\'e}
tel qu'il est {\'e}nonc{\'e} dans \cite{DeMu} 
(voir Th{\'e}or{\`e}me 4.19 et le commentaire qui suit).
Il suffit de v{\'e}rifier la propret{\'e} de $\overline{M^1}$. 

Soit $V$ un anneau de valuation discr{\`e}te de corps de fractions $K$.
Comme $M^1$ est ouvert et dense dans $\overline{M^1}$,  il suffit de v{\'e}rifier que  
tout  morphisme $g^0:\Spec(K)\rightarrow M^1$, s'{\'e}tend en un morphisme 
$g:\Spec(V)\rightarrow \overline{M^1}$. Supposons que 
$g^0$ ne s'{\'e}tend pas d{\'e}j{\`a} en un morphisme  $g:\Spec(V)\rightarrow M^1$.
La VAHB $A/K$ donn{\'e}e par $f^0$ est donc {\`a} mauvaise r{\'e}duction 
(voir Deligne-Pappas \cite{DePa}). 
Quitte {\`a} remplacer $K$ par une extension finie et $V$ par sa normalisation,
on peut supposer que $A/K$ est {\`a} mauvaise r{\'e}duction semi-stable. 
Nous sommes alors en mesure d'appliquer {\`a} $A/K$ la th{\'e}orie de g{\'e}om{\'e}trie rigide de
Raynaud, qui nous fournit (voir la partie \ref{raynaud2}) :

$\relbar$ deux id{\'e}aux $\mathfrak{a}$ et $\mathfrak{b}$, tels 
que $A_{\rig}=(\Gm\otimes \mathfrak{a}^*)_{\rig}/\mathfrak{b}_{\rig}$
et $\mathfrak{c}=\Sym_{\mathfrak{o}}(A,A^t)=\mathfrak{ab}^{-1}$ 
(ceci d{\'e}finit une $R$-pointe).

$\relbar$ une suite exacte 
$0 \rightarrow (\mathfrak{a}/\mathfrak{na})(1) 
\rightarrow A[\mathfrak{n}]\rightarrow 
\mathfrak{n}^{-1}\mathfrak{b}/\mathfrak{b}\rightarrow 0$. La $\mu_{\mathfrak{n}}$-structure de
niveau  $(\mathfrak{o}/\mathfrak{n})(1) \hookrightarrow A[\mathfrak{n}]$
d{\'e}finit alors une $(R,\mathfrak{n})$-composante  $\mathcal{C}$ 
au-dessus de la $R$-pointe
d{\'e}finie pr{\'e}c{\'e}demment ({\`a} laquelle on peut associer 
un id{\'e}al $\mathfrak{b}'\supset \mathfrak{b}$) et une 
racine de l'unit{\'e} $\zeta_{\mathcal{C}}$ d'ordre l'exposant du groupe 
$\mathfrak{b}'/\mathfrak{b}$.

$\relbar$ un  {\'e}l{\'e}ment $\xi^*\in (\mathfrak{ab})^*_+$
(bien d{\'e}fini  modulo l'action de  $\mathfrak{o}_{\mathcal{C}}^{\times}$),
venant de la forme bilin{\'e}aire $\mathfrak{o}$-{\'e}quivariante 
 $\langle,\rangle:\mathfrak{a}\times \mathfrak{b}\rightarrow 
\mathrm{Val(K)}\cong\Z$ \enspace $(\alpha,\beta)\mapsto 
\mathrm{val}(\mathfrak{X}^\alpha(\beta))$.

Un translat{\'e} de  $\xi^*$ par le groupe $\mathfrak{o}_{\mathcal{C},1}^{\times}$
appartient {\`a} un certain c{\^o}ne $\sigma_i^\mathcal{C}\in \Sigma^\mathcal{C}$, 
parmi les c{\^o}nes  choisis dans la d{\'e}monstration du th{\'e}or{\`e}me.
Le morphisme $g^0$ se factorise alors par la carte locale 
$\overline{S}{}_{\sigma_i^\mathcal{C}}^0\times 
\Spec(\Z[\frac{1}{\N(\mathfrak{n})},\zeta_{\mathcal{C}}]) 
\rightarrow  M^1$. Le morphisme compos{\'e} 
$g:\Spec(V) \rightarrow \overline{S}_{\sigma_i^\mathcal{C}}\times
 \Spec(\Z[\frac{1}{\N(\mathfrak{n})},\zeta_{\mathcal{C}}])
\rightarrow \overline{M^1}$ {\'e}tend le morphisme $g^0$. \hfill $\square$

\section{Formes  de Hilbert et compactification minimale arithm{\'e}tiques.}

Nous savons qu'une  forme modulaire de 
Hilbert classique (sur $\C$) est uniquement  d{\'e}termin{\'e}e par son   
$q$-d{\'e}veloppement en une pointe $\mathcal{C}$, que la condition 
d'holomorphie  {\`a} l'infinie est automatiquement  
satisfaite si $d_F>1$ (Principe de Koecher) et qu'il 
n'y a pas de s{\'e}ries d'Eisenstein en poids non-parall{\`e}le.

Le but de cette partie  est de d{\'e}crire, en suivant \cite{rapoport}, les
propri{\'e}t{\'e}s du $q$-d{\'e}veloppement  d'une forme de  Hilbert arithm{\'e}tique.
C'est le point de d{\'e}part de 
la construction de la compactification  minimale arithm{\'e}tique de $M$ 
(voir \cite{chai}). 
 
\subsubsection*{Formes modulaires de Hilbert arithm{\'e}tiques.}\label{fmh}

Pour la  d{\'e}finition de l'espace des formes modulaires de Hilbert,
 nous suivons  le paragraphe  6.8 de \cite{rapoport}, r{\'e}dig{\'e} par P. Deligne.
Consid{\'e}rons le sch{\'e}ma en groupes $\mathcal{T}_1=\Res^{\mathfrak{o}}_{\Z}\Gm$
sur $\Z$, dont la fibr{\'e} g{\'e}n{\'e}rique est le tore $\Res^F_{\Q}\Gm$
de  groupes de caract{\`e}res  $\Z[J_F]$, o{\`u} $J_F$ d{\'e}signe l'ensemble des
plongements de $F$ dans $\R$. On suppose dans cette partie que 
$d_F>1$.

Par d{\'e}finition de la  VAHB universelle, le faisceau
$\underline{\omega}_{\mathcal{A}/M^1}=
f_*\Omega^1_{\mathcal{A}/M^1}$ est un  $\mathfrak{o}$-fibr{\'e} inversible sur
$M^1\times\Spec(\Z[{\frac{1}{\Delta}}])$.
\smallskip

Soit $\kappa\in \Z[J_F]=X(\mathcal{T}_1)$ un poids et soit 
$F''$ un corps de nombres, contenant les valeurs du 
caract{\`e}re $\kappa : F^\times \rightarrow \C^\times$.
On peut prendre, par exemple, $F''=\Q$ et 
poids parall{\`e}le, ou bien $F''=F^{\gal}$
et  poids quelconque.

\smallskip

Soit $\mathfrak{o}''$ l'anneau des entiers de  $F''$. 
Le morphisme de groupes alg{\'e}briques $\kappa : \Res_{\Q}^{F}\Gm \rightarrow
\Res_{\Q}^{F''}\Gm$, se prolonge en un morphisme 
$\Res_{\Z}^{\mathfrak{o}}\Gm \rightarrow \Res_{\Z}^{\mathfrak{o}''}\Gm$, qui 
{\'e}quivaut (par la formule d'adjonction) {\`a} un morphisme de groupes 
alg{\'e}briques sur $\mathfrak{o}''$,
 $\Res_{\Z}^{\mathfrak{o}}\Gm\times\Spec(\mathfrak{o}'')
\rightarrow \Gm\times\Spec(\mathfrak{o}'')$, not{\'e} encore $\kappa$.

\medskip
 
A l'aide de $\kappa$, on peut d{\'e}couper dans $\underline{\omega}$ un 
 fibr{\'e} inversible sur $M^1\times \Spec(\mathfrak{o}''[{\frac{1}{\Delta}}])$, 
not{\'e}   $\underline{\omega}^\kappa$. Soit $\mathfrak{o}'$ 
l'anneau des entiers de $F'=F''(\epsilon^{1/2},\epsilon\in
\mathfrak{o}_{D+}^\times)$. Alors $\underline{\omega}^\kappa$
descend en un fibr{\'e} inversible sur $M\times \Spec(\mathfrak{o}'
[{\frac{1}{\Delta}}])$, not{\'e}  encore $\underline{\omega}^\kappa$ 
(voir la partie 4 de \cite{dimtildg} pour  une pr{\'e}sentation  plus d{\'e}taill{\'e}e). 

Pour tout $\Z[\frac{1}{\Delta}]$-sch{\'e}ma $Y$, on pose 
$Y'=Y\times_{\Spec(\Z[\frac{1}{\Delta}])}
 \Spec(\mathfrak{o}'[\frac{1}{\Delta}])$.

\begin{defin}  Soit $R$ une $\mathfrak{o}'[\frac{1}{\Delta}]$-alg{\`e}bre.
 
Une forme modulaire de Hilbert arithm{\'e}tique 
de poids $\kappa$, de  niveau $\Gamma$ et  {\`a} coefficients dans $R$, est 
une section  globale de $\underline{\omega}^\kappa$ sur
$M\times_{\Spec(\Z[\frac{1}{\Delta}])} \Spec(R)$.
On note 
$G_\kappa(\mathfrak{c},\mathfrak{n};R)^{\mathrm{geom}}:=
\mathrm{H}^0(M\times_{\Spec(\Z[\frac{1}{\Delta}])} \Spec(R)
,\underline{\omega}^{\kappa})$ l'espace  
de ces formes modulaire de Hilbert.
\end{defin}

\begin{rque} 
1) Le faisceau $\underline{\omega}^t$ ($t=\sum_{\tau\in J_F}\tau$)
n'est autre que le faisceau 
$\wedge^{d_F}\underline{\omega}=\det(\underline{\omega})$ 
sur $M$, et $\underline{\omega}^{kt}$ - sa puissance $k$-i{\`e}me.
Les formes modulaires de Hilbert de poids parall{\`e}le $k\geq 1$,
s'{\'e}crivent donc $G_{kt}(\mathfrak{c},\mathfrak{n})^{\mathrm{geom}}=
\mathrm{H}^0(M,(\wedge^{d_F} \underline{\omega})^{\otimes k})$.

\noindent 2) Si $F'\supset F^{\gal }$,  l'action de $\mathfrak{o}$ 
permet de d{\'e}composer  $\underline{\omega}\cong
\mathfrak{o}\otimes\mathcal{O}_{M'} $ en somme directe de fibr{\'e}s inversibles 
$\underline{\omega}^\tau$ sur $M'$, index{\'e}s par les diff{\'e}rents 
plongements $\tau$ de $F$ dans $F^{\gal }$. Si $\kappa=\sum k_\tau\tau$,
on a $\underline{\omega}^\kappa=\otimes_\tau
(\underline{\omega}^\tau)^{\otimes k_\tau}$. 

\end{rque}

 Soit $\overline{f}:\mathfrak{G}\rightarrow \overline{M^1}$ 
le sch{\'e}ma semi-ab{\'e}lien au-dessus d'une compactification 
toro{\"\i}dale $\overline{M^1}$ de $M^1$, comme dans la partie pr{\'e}c{\'e}dente. Le faisceau
 $\underline{\omega}_{\mathfrak{G}/\overline{M^1}}:=
\overline{e}^*\Omega^1_{\mathfrak{G}/\overline{M^1}}$, o{\`u}
$\overline{e}:\overline{M^1}\rightarrow \mathfrak{G}$ d{\'e}signe la
section unit{\'e}, prolonge le faisceau  $\underline{\omega}_{\mathcal{A}/M^1}$.
En outre $\underline{\omega}_{\mathfrak{G}/\overline{M^1}}$ co{\"\i}ncide avec
le faisceau $(\overline{f}_*\Omega^1_{\mathfrak{G}/\overline{M^1}})^{\mathfrak{G}}$
des $\mathfrak{G}$-invariants de 
$\overline{f}_*\Omega^1_{\mathfrak{G}/\overline{M^1}}$.
En passant aux cartes formelles, on voit comme dans \cite{rapoport},
qu'au-dessus de $\Z[{\frac{1}{\Delta}}]$, le faisceau
$\underline{\omega}_{\mathfrak{G}/\overline{M^1}}$ est un 
$\mathfrak{o}$-fibr{\'e} inversible. Le fibr{\'e}
$\underline{\omega}_{\mathfrak{G}/\overline{M^1}}$ descend 
en un $\mathfrak{o}$-fibr{\'e} inversible sur $\overline{M}'$, 
not{\'e} $\underline{\omega}$ (voir la partie 7 de  \cite{dimtildg}).
Pour tout $\kappa\in\Z[J_F]$, on peut ainsi prolonger le fibr{\'e} inversible
$\underline{\omega}^\kappa$ en un fibr{\'e} inversible  sur $\overline{M}'$,
not{\'e} encore $\underline{\omega}^\kappa$.

\smallskip
D'apr{\`e}s la partie \ref{vahb-degenerantes} pour toute
$(R,\mathfrak{n})$-composante uniformis{\'e}e 
$\mathcal{C}$,  tout c{\^o}ne $\sigma\in\Sigma^\mathcal{C}$  et pour toute 
$\mathfrak{o}'[\frac{1}{\Delta},\zeta_{\mathcal{C}}]$-alg{\`e}bre $R$ on a 
$\underline{\omega}|_{S^\wedge_\sigma\times\Spec(R)}\simeq
(\mathfrak{a}\otimes \mathcal{O}_{S^\wedge_\sigma\times
\Spec(R)})$, d'o{\`u} 
\begin{equation}\label{omega}
\underline{\omega}^\kappa|_{S_{\Sigma^\mathcal{C}}^\wedge \times\Spec(R)}\simeq
(\mathfrak{a}\otimes \mathcal{O}_{S_{\Sigma^\mathcal{C}}^\wedge} \otimes
R)^{-\kappa}\!=\!\textstyle
(\mathfrak{a}\otimes \mathfrak{o}'[\frac{1}{\Delta}])^{-\kappa}
\!\underset{\mathfrak{o}'[\frac{1}{\Delta}]}{\otimes} \!
(\mathfrak{o}\otimes \mathcal{O}_{S_{\Sigma^\mathcal{C}}^\wedge} \otimes
R)^{-\kappa} 
\end{equation}

Par cons{\'e}quent 
$\mathrm{H}^0((S_{\Sigma^\mathcal{C}}^\wedge \times\Spec(R))/
\mathfrak{o}_{\mathcal{C}}^\times,\underline{\omega}^\kappa)=\mathfrak{a}^{(\kappa)}
\otimes_{\mathfrak{o}'[\frac{1}{\Delta}]} R_\mathcal{C}^{(\kappa)}(R)$, o{\`u} 
 $\mathfrak{a}^{(\kappa)}=
(\mathfrak{a}\otimes \mathfrak{o}'[\frac{1}{\Delta}])^{-\kappa}$ 
est un $\mathfrak{o}'[\frac{1}{\Delta}]$-module inversible et 
$$R_\mathcal{C}^{(\kappa)}(R):=
\left\{ \sum_{\xi\in X_+\cup\{0\}}\!\! a_\xi q^\xi\Big|
a_\xi\in R, \enspace a_{u^2\epsilon \xi}= \epsilon^{\kappa/2} u^\kappa 
\zeta_{\mathcal{C}}^{n\Tr_{F\!/\!\Q}( \xi u\xi_{u,\epsilon}^*)}a_\xi
, \enspace \forall (u,\epsilon)\in\mathfrak{o}_{\mathcal{C}}^\times\right\}.$$

Notons que $\xi u\xi_{u,\epsilon}^*$ est un {\'e}l{\'e}ment de
 $\mathfrak{b'b}^{-1}\mathfrak{d}^{-1}$, 
bien d{\'e}fini modulo $\mathfrak{d}^{-1}$, et donc
$n\Tr_{F\!/\!\Q}( \xi u\xi_{u,\epsilon}^*)\in \Z/n\Z$ 
(on rappelle que $n\Z=\Z\cap \mathfrak{bb'}^{-1}$
et $n=\mathrm{ord}(\zeta_{\mathcal{C}})$).

On a $R_\mathcal{C}^{(\kappa)}(R)=
 \mathrm{H}^0\big(S_{\Sigma^\mathcal{C}}^\wedge\times \Spec(R),
(\mathfrak{o}\otimes 
\mathcal{O}_{S_{\Sigma^\mathcal{C}}^\wedge}\otimes R)^{-\kappa}
\big)^{\mathfrak{o}_{\mathcal{C}}^\times}$.

Le diagramme suivant montre comment l'anneau 
$R_\mathcal{C}^{(\kappa)}(R)$ 
se situe par rapport aux diff{\'e}rents anneaux d{\'e}j{\`a} consid{\'e}r{\'e}s dans 
la partie \ref{vahb-degenerantes} :

$$\xymatrix{& R[q^\xi]_{\xi\in X_+\cup\{0\}}\ar@{^{(}->}[r]\ar@{^{(}->}[d]
& R_\sigma\otimes R=R[q^\xi]_{\xi\in X\cap\check{\sigma}} \ar@{^{(}->}[d]\\
R_\mathcal{C}^{(\kappa)}(R)\ar@{^{(}->}[r] &
R[[q^\xi]]_{\xi\in X_+\cup\{0\}} \ar@{^{(}->}[r] &
R_\sigma^\wedge\otimes R }$$

\vspace{-.4cm} 
\subsubsection*{Principe de Koecher.}

\begin{theo}{\rm (Principe de Koecher \cite{rapoport} 4.9)}
Soit $\overline{M}$ une compactification toro{\"\i}dale de $M$. Alors 
$$ \mathrm{H}^0\big(M \times \Spec(R), \underline{\omega}^\kappa \big)=
 \mathrm{H}^0\big(\overline{M}\times \Spec(R), \underline{\omega}^\kappa \big)$$
\end{theo}

\noindent {\bf D{\'e}monstration : } Le probl{\`e}me est local et 
il suffit de le v{\'e}rifier apr{\`e}s compl{\'e}tion, le long d'une 
 $(R,\mathfrak{n})$-composante $\mathcal{C}$. 

D'apr{\`e}s la trivialisation  (\ref{omega}) 
du fibr{\'e} inversible  $\underline{\omega}^\kappa$, 
 il  s'agit de voir que les sections globales m{\'e}romorphes 
du faisceau  $(\mathfrak{o}\otimes 
\mathcal{O}_{S_{\Sigma^\mathcal{C}}^\wedge}\otimes R)^{-\kappa}$
sur $S_{\Sigma^\mathcal{C}}^\wedge\times \Spec(R)$, qui sont 
$\mathfrak{o}_{\mathcal{C}}^\times$-invariantes, 
appartiennent  {\`a}  $R_{\mathcal{C}}^{(\kappa)}(R)$.

 Soit $f=\sum_{\xi\in X}a_\xi q^\xi \in \mathrm{H}^0_{\mathrm{mer}}
((S_{\Sigma^\mathcal{C}}^\wedge\times 
\Spec(R))/\mathfrak{o}_{\mathcal{C}}^\times, (\mathfrak{o}\otimes 
\mathcal{O}_{S_{\Sigma^\mathcal{C}}^\wedge}\otimes R)^{-\kappa})$
une telle section. 
Supposons que $a_{\xi_0}\neq 0$ avec $\xi_0$ non-totalement positif. 
Il existe donc $\xi_0^*\in X^*_{\R+}$ avec $\Tr_{F\!/\!\Q}(\xi_0 \xi_0^*)$ 
strictement n{\'e}gatif. Comme $d_F>1$, 
on peut choisir des unit{\'e}s  $u \in \mathfrak{o}_{\mathcal{C},1}^\times$ 
de mani{\`e}re {\`a} rendre la quantit{\'e} $\Tr_{F\!/\!\Q}(u^2\xi_0 \xi_0^*)$ 
arbitrairement  proche de $-\infty$.   Soit $\sigma$ un c{\^o}ne poly{\'e}dral 
de $\Sigma^{\mathcal{C}}$ contenant $\xi_0^*$. Par d{\'e}finition de 
$S_{\sigma}^\wedge$, on  voit que $f$ n'est pas
m{\'e}romorphe sur $S_{\sigma}^\wedge$. Contradiction.
Donc $f\in R_{\mathcal{C}}^{(\kappa)}(R)$. \hfill $\square$

\subsubsection*{$q$-d{\'e}veloppement.}

Le paragraphe  pr{\'e}c{\'e}dent montre que l'on peut associer  {\`a} une 
$(R,\mathfrak{n})$-composante  uniformis{\'e}e $\mathcal{C}$ et {\`a} une 
forme modulaire de Hilbert $f$ de poids $\kappa$,  
niveau $\Gamma$, et {\`a} coefficients dans une 
$\mathfrak{o}'[\frac{1}{\Delta},\zeta_\mathcal{C}]$-alg{\`e}bre  $R$, 
un {\'e}l{\'e}ment :

$$f_\mathcal{C}\in  \mathfrak{a}^{(\kappa)} 
\otimes_{\mathfrak{o}'[\frac{1}{\Delta}]}
R_\mathcal{C}^{(\kappa)}(R).$$

\begin{defin}
L'{\'e}l{\'e}ment $f_\mathcal{C}$  est appel{\'e} le $q$-d{\'e}veloppement 
de la forme $f$ 
en la $(R,\mathfrak{n})$-composante 
uniformis{\'e}e  $\mathcal{C}$.
On note $\ev_{\mathcal{C},\kappa}$ l'application $f\mapsto f_\mathcal{C}$.
\end{defin}

 Le principe du $q$-d{\'e}veloppement s'{\'e}nonce alors :

\begin{prop} Soient $\kappa$ un poids, 
$\mathcal{C}$ une $(R,\mathfrak{n})$-composante 
uniformis{\'e}e et  $R$ une 
$\mathfrak{o}'[\frac{1}{\Delta},\zeta_{\mathcal{C}}]$-alg{\`e}bre
(contenant les valeurs de $\kappa$). Alors 

\noindent $\mathrm{(i)}$ l'application $\ev_{\mathcal{C},\kappa}$ est
injective,

\noindent$\mathrm{(ii)}$ pour  toute  $R$-alg{\`e}bre  $R'$  et
 $f\in G_\kappa(\mathfrak{c},\mathfrak{n};R')$, 
si   $\ev_{\mathcal{C},\kappa}(f)\in 
\mathfrak{a}^{(\kappa)}\otimes_{\mathfrak{o}'[\frac{1}{\Delta}]}
R_\mathcal{C}^{(\kappa)}(R)$,
alors $f\in G_\kappa(\mathfrak{c},\mathfrak{n};R)$,

\noindent$\mathrm{(iii)}$ s'il existe  $f\in G_\kappa(\mathfrak{c},\mathfrak{n};R)$
tel que le terme constant de  
$\ev_{\mathcal{C},\kappa}(f)$ ne soit  pas nul,  alors 
$\epsilon^{\kappa/2} u^\kappa -1$ 
est un diviseur de z{\'e}ro dans $R$, pour tout 
$(u,\epsilon) \in \mathfrak{o}_{\mathcal{C}}^\times$.

\end{prop}

Le cas de l'anneau nul $R=0$ redonne une formulation classique du principe.
Pour d{\'e}monstration du  (i) et du (ii)  
voir la partie 7 de \cite{dimtildg}. Le (iii) est clair. \hfill $\square$

\subsubsection*{Compactification minimale.}

La compactification minimale est la contrepartie arithm{\'e}tique de 
la compactification de Satake sur $\C$. Contrairement au cas complexe,
dans le cas arithm{\'e}tique, la   construction  de la compactification minimale
utilise les compactifications toro{\"\i}dales.  Voici l'analogue en niveau
 $\Gamma$ de l'{\'e}nonc{\'e} donn{\'e}   par C.-L. Chai dans \cite{chai}.

\begin{theo} ${}$

\noindent $(\mathrm{i})$ Il existe $k_0\in \n^*$ tel que 
le faisceau $\underline{\omega}_{\mathcal{A}/M^1}^{k_0t}$, 
soit engendr{\'e} par ses sections globales sur $\overline{M^1}$.

\noindent $(\mathrm{ii})$ Le morphisme canonique 
$\pi :\overline{M^1}\rightarrow M^{1*}:=\Proj_{\Z[\frac{1}{\N(\mathfrak{n})}]}
\left(\oplus_{k\geq 0}  \mathrm{H}^0(\overline{M^1},
\underline{\omega}_{\mathcal{A}/M^1}^{kt})\right), $
est surjectif. Le $\Z[\frac{1}{\N(\mathfrak{n})}]$-sch{\'e}ma $M^{1*}$
est ind{\'e}pendant du choix de  $\Sigma$ (on rappelle que $\overline{M^1}=M^1_\Sigma$).

\noindent $(\mathrm{iii})$ L'anneau gradu{\'e} $\oplus_{k\geq 0}
 \mathrm{H}^0(\overline{M^1},\underline{\omega}_{\mathcal{A}/M^1}^{kt})$ est de type fini sur 
$\Z[\frac{1}{\N(\mathfrak{n})}]$ et $M^{1*}$ est un $\Z[\frac{1}{\N(\mathfrak{n})}]$-sch{\'e}ma projectif,
normal, de type fini.  Le groupe  $\mathfrak{o}_{D+}^\times/
(\mathfrak{o}_{D+}^\times\cap\mathfrak{o}_{\mathfrak{n}}^{\times 2})$
du rev{\^e}tement fini {\'e}tale $\overline{M^1}\rightarrow \overline{M}$ 
agit sur $M^{1*}$ et le quotient est un 
$\Z[\frac{1}{\N(\mathfrak{n})}]$-sch{\'e}ma 
projectif, normal, de type fini $M^*$, muni d'un morphisme 
surjectif $\pi :\overline{M}\rightarrow M^{*}$. 

\noindent $(\mathrm{iv})$ $\pi|_{M}$ induit   un isomorphisme
sur un ouvert dense de $M^{*}$, not{\'e} encore $M$. 
$M^{*} \backslash M$ est fini et {\'e}tale sur 
$\Z[\frac{1}{\N(\mathfrak{n})}]$ et en fait
isomorphe {\`a} : $$ \coprod_{(R,\mathfrak{n})\mathrm{-composantes}/\sim}
 \textstyle
\Spec(\Z[\frac{1}{\N(\mathfrak{n})},\zeta_{\mathcal{C}}]^{H_{\mathcal{C}}}).$$
Les composantes connexes de $M^{*} \backslash M$ sont appel{\'e}es les
pointes de $M$. Cependant celles-ci ne sont des points ferm{\'e}s que pour 
les $(R,\mathfrak{n})$-composantes non-ramifi{\'e}es.

\noindent $(\mathrm{v})$ L'image r{\'e}ciproque $\pi^{-1}(\mathcal{C})$
de chaque pointe $\mathcal{C}$ de $M$, est une composante connexe
de $\overline{M} \backslash M$.  La compl{\'e}tion formelle de 
$\overline{M}$ le long de l'image r{\'e}ciproque d'une 
composante $\pi^{-1}(\mathcal{C})$, est canoniquement isomorphe {\`a} 
$(S_{\Sigma^{\mathcal{C}}}^{\wedge}/\mathfrak{o}_{\mathcal{C}}^{\times})
\times  \Spec(\Z[\frac{1}{\N(\mathfrak{n})},
\zeta_{\mathcal{C}}]^{H_{\mathcal{C}}})$.
En particulier, la compl{\'e}tion formelle de 
$\overline{M}$ le long de l'image r{\'e}ciproque $\pi^{-1}(\mathcal{C})$ d'une   
$(R,\mathfrak{n})$-composante non-ramifi{\'e}e $\mathcal{C}$, est canoniquement
isomorphe {\`a} 
$$(S_{\Sigma^{\mathcal{C}}}^{\wedge}/
(\mathfrak{o}_{\mathfrak{n}}^{\times}\times\mathfrak{o}_{D+}^\times))
{ }  \times \textstyle \Spec(\Z[\frac{1}{\N(\mathfrak{n})}]).$$

\noindent $(\mathrm{vi})$ Pour tout $\kappa\in\Z[J_F]$ le faisceau 
$\underline{\omega}^{\kappa}$ s'{\'e}tend en un faisceau 
inversible sur $M^{*}$  si et seulement si $\kappa$ est parall{\`e}le.
\end{theo}

\noindent {\bf D{\'e}monstration : } Nous suivons  la m{\'e}thode de C.-L. Chai \cite{chai}.
D'apr{\`e}s \cite{MB} Chap.IX Thm.2.1 (voir aussi \cite{FaCh} Chap.V Prop.2.1), 
il existe $k_0\geq 1$ tel que le faisceau 
inversible $\underline{\omega}_{\mathcal{A}/M^1}^{k_0t}$ 
soit engendr{\'e} par ses sections globales sur $\overline{M^1}$.
Ceci nous fournit un morphisme 
$$\overline{M^1}\rightarrow \Proj_{\Z[\frac{1}{\N(\mathfrak{n})}]}\left(\Sym^{\bullet}
 \mathrm{H}^0(\overline{M^1},\underline{\omega}_{\mathcal{A}/M^1}^{k_0t})\right).$$

Soit $B^{\bullet}$ la normalisation de $\Sym^{\bullet}
 \mathrm{H}^0(\overline{M^1},\underline{\omega}_{\mathcal{A}/M^1}^{kk_0t})$ dans  $\underset{k\geq 0}{\oplus}
 \mathrm{H}^0(\overline{M^1},\underline{\omega}_{\mathcal{A}/M^1}^{kk_0t})$. Le morphisme
associ{\'e} $\pi:\overline{M^1}\rightarrow \Proj_{\Z[\frac{1}{\N(\mathfrak{n})}]}
(B^{\bullet})$ est birationnel, surjectif et satisfait $\pi^*\mathcal{O}(1)=\underline{\omega}_{\mathcal{A}/M^1}^{k_0t}$.
Le Th{\'e}or{\`e}me de Connexit{\'e} de Zariski implique alors que les fibres de
$\pi$ sont connexes. D'apr{\`e}s  \cite{FaCh} Chap.V Prop.2.2, la partie 
ab{\'e}lienne est constante dans chaque fibre g{\'e}om{\'e}trique de $\pi$, et par 
cons{\'e}quent les fibres g{\'e}om{\'e}triques de $\pi$ sont 

\noindent $\relbar$ soit des points g{\'e}om{\'e}triques de $M^1$,

\noindent $\relbar$ soit des composantes g{\'e}om{\'e}triques connexes de 
$\overline{M^1}\bs M^1$.

\smallskip

Comme pour tout $k\geq 1$,  $\pi^*\mathcal{O}(k)=\underline{\omega}_{\mathcal{A}/M^1}^{k_0kt}$ et 
$\underline{\omega}_{\mathcal{A}/M^1}^{k_0t}$ est engendr{\'e} par ses sections globales sur $\overline{M^1}$, 
on obtient $ \mathrm{H}^0(\overline{M^1},\underline{\omega}_{\mathcal{A}/M^1}^{kk_0t})=
 \mathrm{H}^0(\Proj(B^{\bullet}),\mathcal{O}(k))$. Par cons{\'e}quent 
$B^{\bullet}=\underset{k\geq 0}{\oplus} \mathrm{H}^0(\overline{M^1},\underline{\omega}_{\mathcal{A}/M^1}^{kk_0t})$ et 
c'est une $\Z[\frac{1}{\N(\mathfrak{n})}]$-alg{\`e}bre de type fini. 
Or, l'alg{\`e}bre  $\underset{k\geq 0}{\oplus}  \mathrm{H}^0(\overline{M^1},\underline{\omega}_{\mathcal{A}/M^1}^{kt})$
est enti{\`e}re sur $\underset{k\geq 0}{\oplus} 
\mathrm{H}^0(\overline{M^1},\underline{\omega}_{\mathcal{A}/M^1}^{kk_0t})$,
engendr{\'e}e par les {\'e}l{\'e}ments de degr{\'e} plus petit que $k_0$. Il en r{\'e}sulte que
$\underset{k\geq 0}{\oplus}  \mathrm{H}^0(\overline{M^1},\underline{\omega}_{\mathcal{A}/M^1}^{kt})$ est de type fini  sur $\Z[\frac{1}{\N(\mathfrak{n})}]$, 
et que  $M^{1*}:=\Proj(\underset{k\geq 0}{\oplus}  
\mathrm{H}^0(\overline{M^1},\underline{\omega}_{\mathcal{A}/M^1}^{kt}))
=\Proj(B^{\bullet})$. Par le principe de
Koecher, le sch{\'e}ma $M^{1*}$ est ind{\'e}pendant du choix particulier 
de la  compactification  toro{\"\i}dale $\overline{M^1}$ de $M^1$. 
Le groupe  $\mathfrak{o}_{D+}^\times/
(\mathfrak{o}_{D+}^\times\cap\mathfrak{o}_{\mathfrak{n}}^{\times 2})$
agit sur $M^{1*}$ et on d{\'e}finit $M^*$ comme le quotient. 
Notons qu'en g{\'e}n{\'e}ral $M^{1*}\rightarrow M^*$ n'est pas {\'e}tale, car 
les pointes  peuvent avoir des stabilisateurs non-triviaux. 

On a donc (i),(ii),(iii) et la premi{\`e}re partie de (iv).  Le calcul de 
la compl{\'e}tion formelle de $\overline{M}$, le long de l'image r{\'e}ciproque 
d'une composante connexe de $M^{*}\bs M$ d{\'e}coule du Th{\'e}or{\`e}me des
 Fonctions Formelles de Grothendieck. 

Enfin, examinons {\`a} quelle condition   
$\underline{\omega}_{\mathfrak{G}\!/\!\overline{M^1}}^\kappa$ descend en un 
fibr{\'e} inversible sur $M^{1*}$.
Comme $(\pi_*\underline{\omega}_{\mathfrak{G}\!/\!\overline{M^1}}^\kappa)
|_{M^1}=
\underline{\omega}_{\mathcal{A}/M^1}^\kappa$ et 
$\mathrm{codim}(M^{1*}\bs M^1)\geq 2$, 
le faisceau $\pi_*\underline{\omega}_{\mathfrak{G}\!/\!\overline{M^1}}^\kappa$ est coh{\'e}rent.
 Il est inversible si et seulement si 
$\underline{\omega}_{\mathfrak{G}\!/\!\overline{M^1}}^\kappa$
peut {\^e}tre trivialis{\'e} sur  $S_{\Sigma^{\mathcal{C}}}^{\wedge}/
\mathfrak{o}_{\mathcal{C},1}^{\times}\times  \Spec(R)$.
D'apr{\`e}s (\ref{omega}) le pull-back de $\underline{\omega}_{\mathfrak{G}\!/
\!\overline{M^1}}^\kappa$  {\`a} $S_{\Sigma^{\mathcal{C}}}^{\wedge}\times \Spec(R)$
est canoniquement trivial et  $\mathfrak{o}_{\mathcal{C},1}^{\times}$ agit sur 
ce pull-back {\`a} travers $\kappa$, d'o{\`u}  (vi). 
\hfill $\square$

\vspace{-.4cm} 
\subsubsection*{Exemples de $q$-d{\'e}veloppement en une pointe ramifi{\'e}e.}

Nous nous proposons de d{\'e}crire explicitement dans le cas
particulier  de l'exemple \ref{Rn-exemple}(ii)(iii) 
les anneaux  $R_\mathcal{C}^{(\kappa)}(R)$
contenant les $q$-d{\'e}veloppements 
des formes modulaires de Hilbert de poids $\kappa$ et  niveau $\Gamma$.
 Rappelons que $\mathfrak{o}'$ d{\'e}signe les entiers d'un corps de nombres 
contenant les valeurs du caract{\`e}re $\kappa$. On suppose que
$\mathrm{Cl}_F=\{1\}$.

\medskip
Pla{\c c}ons nous dans  le cas (ii). Le bord $M^{1*}\bs M^1$ s'{\'e}crit alors 
$$\mathop{\coprod_{(R,\mathfrak{n})\mathrm{-comp.}}}_{\text{non-ramifi{\'e}s}/\sim}  \textstyle
\hspace{-5mm}\Spec\left(\Z\left[\frac{1}{\N(\mathfrak{n})}\right]\right)\hspace{-4mm}\displaystyle
\mathop{\coprod_{(R,\mathfrak{n})\mathrm{-comp.}}}_{\text{peu ramifi{\'e}s}/\sim}  \textstyle
\hspace{-5mm}\Spec\left(\Z\left[\frac{1}{\N(\mathfrak{n})},\zeta_p\right]\right)\hspace{-4mm}\displaystyle
\mathop{\coprod_{(R,\mathfrak{n})\mathrm{-comp.}}}_{\text{tr{\`e}s ramifi{\'e}s}/\sim}  \textstyle
\hspace{-7mm}\Spec\left(\Z\left[\frac{1}{\N(\mathfrak{n})},\zeta_{p^2}\right]^{\mathfrak{o}^\times\!/\!\mathfrak{o}^\times_{\mathfrak{p}^2}}\right)$$

\noindent $\relbar$  Si la pointe $\mathcal{C}$ est non-ramifi{\'e}e, 
pour toute $\mathfrak{o}'[\frac{1}{\Delta}]$-alg{\`e}bre 
$R$, on a $$R_\mathcal{C}^{(\kappa)}(R)=\left\{ \sum_{\xi\in \mathfrak{o}_+} a_\xi q^\xi |
a_\xi\in R, \enspace a_{u^2\xi}= u^\kappa  a_\xi
, \enspace \forall u\in\mathfrak{o}_{\mathfrak{p}^2}^\times\right\}.$$

\noindent $\relbar$  Si la pointe $\mathcal{C}$ est peu ramifi{\'e}e, 
pour toute $\mathfrak{o}'[\frac{1}{\Delta},\zeta_p]$-alg{\`e}bre 
$R$, on a $$R_\mathcal{C}^{(\kappa)}(R)=\left\{ \sum_{\xi\in \mathfrak{p}^{-1}_+} a_\xi q^\xi |
a_\xi\in R, \enspace a_{u^2\xi}= u^\kappa 
\zeta_{p}^{p\Tr_{F\!/\!\Q}( \xi u\xi_u^*)}a_\xi
, \enspace \forall u\in\mathfrak{o}_{\mathfrak{p}}^\times\right\}.$$

\noindent $\relbar$  Si la pointe $\mathcal{C}$ est tr{\`e}s ramifi{\'e}e, 
pour toute $\mathfrak{o}'[\frac{1}{\Delta},\zeta_{p^2}]$-alg{\`e}bre 
$R$, on a 
$$R_\mathcal{C}^{(\kappa)}(R)= 
\left\{ \sum_{\xi\in \mathfrak{p}^{-2}_+} a_\xi q^\xi |
a_\xi\in R, \enspace a_{u^2\xi}=
 u^\kappa \zeta_{p^2}^{p^2\Tr_{F\!/\!\Q}( \xi u\xi_u^*)} a_\xi , 
\enspace \forall u\in\mathfrak{o}_{\mathfrak{p}^2}^\times\right\}.$$

En fait, d'apr{\`e}s la d{\'e}monstration de la Prop.\ref{deuxideaux}(iv), 
on a  $\xi_u^*\in \mathfrak{p}^2\mathfrak{d}^{-1}$ et donc 
$\Tr_{F\!/\!\Q}( \xi u\xi_u^*)\in \Z$ (alors qu'{\`a} priori 
il appartient juste {\`a} $p^{-2}\Z$!). On en d{\'e}duit que 

$$R_\mathcal{C}^{(\kappa)}(R)= 
\left\{ \sum_{\xi\in \mathfrak{p}^{-2}_+} a_\xi q^\xi |
a_\xi\in R, \enspace a_{u^2\xi}= u^\kappa  a_\xi , 
\enspace \forall u\in\mathfrak{o}_{\mathfrak{p}^2}^\times\right\},$$
ce qui est compatible avec le fait que $\zeta_{p^2}$
n'appartient pas au corps de d{\'e}finition de la pointe $\mathcal{C}$, 
qui est $\Q(\zeta_{p^2})^{\mathfrak{o}^\times\!/\!\mathfrak{o}^\times_{\mathfrak{p}^2}}$ (notons que $-1\in\mathfrak{o}^\times\!/\!\mathfrak{o}^\times_{\mathfrak{p}^2}$).

\bigskip
Pla{\c c}ons nous dans  le cas (iii). Le bord $M^{1*}\bs M^1$ s'{\'e}crit alors 
$$
\mathop{\coprod_{(R,\mathfrak{n})\mathrm{-comp.}}}_{\text{non-ramifi{\'e}s}/\sim}
\textstyle\Spec\left(\Z\left[\frac{1}{\N(\mathfrak{n})}\right]\right)
\displaystyle
\mathop{\coprod_{(R,\mathfrak{n})\mathrm{-comp.}}}_{\text{ramifi{\'e}s}/\sim} 
\textstyle\Spec\left(\Z\left[\frac{1}{\N(\mathfrak{n})},\zeta_p\right]^{\mathfrak{o}_{\overline{\mathcal{C}}}^\times/
\mathfrak{o}^\times_{\mathfrak{p}}}\right).$$

\noindent $\relbar$  Si la pointe $\mathcal{C}$ est non-ramifi{\'e}e, 
pour toute $\mathfrak{o}'[\frac{1}{\Delta}]$-alg{\`e}bre 
$R$, on a $$R_\mathcal{C}^{(\kappa)}(R)=\left\{ \sum_{\xi\in \mathfrak{o}_+} a_\xi q^\xi |
a_\xi\in R, \enspace a_{u^2\xi}= u^\kappa  a_\xi
, \enspace \forall u\in\mathfrak{o}_{\mathfrak{p}}^\times\right\}.$$

\noindent $\relbar$  Si la pointe $\mathcal{C}$ est  ramifi{\'e}e, 
pour toute $\mathfrak{o}'[\frac{1}{\Delta},\zeta_{p}]$-alg{\`e}bre 
$R$, on a 
$$R_\mathcal{C}^{(\kappa)}(R)=\left\{ \sum_{\xi\in \mathfrak{p}^{-1}_+} a_\xi q^\xi |
a_\xi\in R, \enspace a_{u^2\xi}= u^\kappa \zeta_{p}^{p\Tr_{F\!/\!\Q}( \xi u\xi_u^*)} 
a_\xi , \enspace \forall u\in\mathfrak{o}_{\mathfrak{p}}^\times\right\}.$$

En fait, d'apr{\`e}s la d{\'e}monstration de la Prop.\ref{deuxideaux}(iv), 
on a $\xi_u^*\in \mathfrak{pd}^{-1}$  et donc 
$\Tr_{F\!/\!\Q}( \xi u\xi_u^*)\in \Z$ (alors qu'{\`a} priori 
il appartient juste {\`a} $p^{-1}\Z$!). On en d{\'e}duit que 

$$R_\mathcal{C}^{(\kappa)}(R)=\left\{ \sum_{\xi\in \mathfrak{p}^{-1}_+} a_\xi q^\xi |
a_\xi\in R, \enspace a_{u^2\xi}= u^\kappa  a_\xi , \enspace \forall u\in\mathfrak{o}_{\mathfrak{p}}^\times\right\},$$
ce qui est compatible avec le fait que $\zeta_{p}$
n'appartient pas au corps de d{\'e}finition de la pointe $\mathcal{C}$, 
qui est $\Q(\zeta_{p})^{\mathfrak{o}_{\overline{\mathcal{C}}}^\times/
\mathfrak{o}^\times_{\mathfrak{p}}}$ (notons que $-1\in
\mathfrak{o}_{\overline{\mathcal{C}}}^\times/
\mathfrak{o}^\times_{\mathfrak{p}}$).

\bibliographystyle{siam}

\bigskip

Universit{\'e} Paris 7 

UFR de Math{\'e}matiques,  Case 7012

2 place Jussieu

75251 PARIS 

FRANCE

\medskip
 Email : \texttt{dimitrov@math.jussieu.fr}

\end{document}